\documentclass[10pt%,twoside
]{article}
\usepackage{latexsym}
\usepackage{amsthm}
\usepackage{amsfonts}
\usepackage{amssymb}
\usepackage{color}

\theoremstyle{plain}
  \newtheorem{theorem}{Theorem}[section]
  
  \newtheorem{proposition}{Proposition}[section]
  \newtheorem{lemma}{Lemma}[section]
  \newtheorem{corollary}{Corollary}[section]

\theoremstyle{remark}
  \newtheorem{remark}{Remark}[section]

\theoremstyle{definition}
  \newtheorem{definition}{Definition}[section]
  \newtheorem{notation}{Notation}[section]

\makeatletter

\@addtoreset{equation}{section}
\makeatother

\begin{document}

\title{
Energy transfer model for the derivative nonlinear Schr\"{o}dinger equations on the torus}
\author{Hideo Takaoka\thanks{This work was supported by JSPS KAKENHI Grant Number 10322794.}\\
Department of Mathematics, Hokkaido University\\
Sapporo, 060-0810, Japan\\
takaoka@math.sci.hokudai.ac.jp}

\pagestyle{myheadings}
\markboth{Hideo Takaoka}{Derivative nonlinear Schr\"odinger equation on torus}

\date{}
      
\maketitle

\begin{abstract}
We consider the nonlinear derivative Schr\"odinger equation with a quintic nonlinearity, on the one dimensional torus.
We exhibit that the nonlinear dynamic properties consist of four frequency modes initially excited, whose frequencies include the resonant clusters and phase matched resonant interactions of nonlinearities.
This phenomena arrests energy transfers between low and high modes, which are quantified by a growth in the Sobolev norm. 
\end{abstract}

\noindent
{\it $2010$ Mathematics Subject Classification{\rm :}} 
Primary  35Q55;  Secondary 42B37.

\noindent
{\it Key words and phrases{\rm :}}
nonlinear Schr\"odinger equation; energy transfer
%\end{abstract}

\section{Introduction}\label{intro}
\indent

In this paper we consider the derivative nonlinear Schr\"odinger equation with quintic nonlinearity:
\begin{eqnarray}\label{dnls}
\left\{
\begin{array}{l}
i\partial u+\partial_x^2u=-i\lambda u^2\partial_x\overline{u}+\mu|u|^4u,\quad (t,x)\in \mathbb{R}\times \mathbb{T},\\
u(0,x)=u_0(x),\quad x\in \mathbb{T},
\end{array}
\right.
\end{eqnarray}
where $u=u(t,x):\mathbb{R}\times \mathbb{T}\to \mathbb{C},~\lambda,~\mu \in \mathbb{R}$, and $\mathbb{T}=\mathbb{R}/2\pi \mathbb{Z}$ is the torus.
Our aim is to prove that there exist solutions of (\ref{dnls}) which initially oscillate, and illustrate the energy transfer from one coupled oscillator to another.

The derivative nonlinear Schr\"odinger equation in the completely integrable form is written as
\begin{eqnarray}\label{dnls-ori}
i\partial u+\partial_x^2u=i\lambda\partial_x(|u|^2u).
\end{eqnarray}
As in \cite{herr}, by using the gauge transform
\begin{eqnarray*}%\label{gauge}
v(t,x)=u(t,x)e^{-i{\cal G}[u]},
\end{eqnarray*}
%where
$$
{\cal G}[u](t,x)=\frac{\lambda}{2\pi}\int_0^{2\pi}\!\int_{\theta}^x\left(|u(t,y)|^2-\frac{1}{2\pi}\|u(t)\|_{L_x^2}^2\right)\,dyd\theta,
$$
we derive a equivalent equation of (\ref{dnls-ori}) for $v$ as follows:
\begin{eqnarray}
i\partial_tv+\partial_x^2v -i\frac{\lambda}{\pi}\|v(t)\|_{L_x^2}^2v& = & -i\lambda v^2\partial_x\overline{v}-\frac{\lambda^2}{2}|v|^4v\nonumber\\
 & & -\frac{\lambda}{\pi}\left(\int_0^{2\pi}\mathrm{Im}(v\partial_x\overline{v})(t,\theta)\,d\theta\right)v\nonumber\\
 & & +\frac{\lambda^2}{2\pi}\|v(t)\|_{L_x^2}^2|v|^2v\nonumber\\
& & +\frac{\lambda^2}{4\pi}\left(\|v(t)\|_{L_x^4}^4-\frac{1}{\pi}\|v(t)\|_{L_x^2}^4\right)v.\label{gauged-dnls-ori}
\end{eqnarray}
The equations (\ref{dnls-ori}) and (\ref{gauged-dnls-ori}) are completely integrable Hamiltonian equations, having infinitely many conservation laws \cite{cc}.

The well-posedness results for the equation (\ref{dnls-ori}), in particular considering solutions in the Sobolev space $H^s$, were studied in many works (e.g., \cite{herr,ta2,su}).
If $u_0\in H^s$ with $s\ge 1/2$,  local well-posedness was established by Herr \cite{herr}.
This combined with the conservation laws, can often yield the corresponding global well-posedness result for data in $H^s$ for $s\ge 1$, which is small in $L^2$.
Extension of the global well-posedness with low regularity data was obtained in the case $s>1/2$ by Su Win \cite{su}, by using the I-method developed in \cite{CKSTT0,ckstt1,CKSTT}.
A further extension was given by \cite{GH,nors,ta2}. 
For instance, in \cite{GH} Gr\"unrock and Herr proved the local well-posedness result for initial data in Fourier-Lebesgue space $\widehat{H}^{s}_r$ (\ref{dnls}), in the range $s\ge 1/2,~4/3<r<2$, where the space $\widehat{H}^{s}_r$ is defined by the norm
$$
\|f\|_{\widehat{H}^s_r}=\left\|\langle \xi\rangle^s\widehat{f}\right\|_{\ell_{\xi}^{r'}},\quad \langle\xi\rangle=(1+\xi^2)^{1/2},\quad \frac{1}{r}+\frac{1}{r'}=1.
$$
Later on in \cite{nors} Nahmod, Oh, Rey-Bellt and Staffilani established the almost sure global well-posedness of (\ref{dnls-ori}) with random initial data in $\widehat{H}^{s}_r$ with $s>1/2$ and certain $4/3<r<2$.
Below $s=1/2$, in \cite{ta2} Takaoka obtained the existence and time continuity of solution in $H^s$ with some $s<1/2$.
We recall that the initial value problem of the equation (\ref{dnls-ori}) is ill-posed in $H^s$ provided $s<1/2$ (see \cite{bialin}).
In the case of the line case setting, analogous results to the well-posedness of the equation (\ref{dnls-ori}) were taken from the references \cite{ckstt1,CKSTT,ho1,ho2,my,oz,ta}.
 
The above technology can also be applied to the equation in (\ref{dnls}) to obtain the local and global well-posedness in  $H^s$.
Indeed, the initial value problem (\ref{dnls}) is locally well-posed in $H^s$ for all $s\ge 1/2$.
In order to capture the global solution in the energy class $s\ge 1$, we need a priori estimate. 
Formally, solutions of (\ref{dnls}) satisfy the conservations laws, that is, the following three quantities are constant in time:
\begin{eqnarray}\label{conservation}
\begin{array}{ll}
M[u](t)=\frac{1}{2\pi}\int_0^{2\pi}m[u](t,x)\,dx=M[u](0) &  \mbox{($L^2$-norm)},\\
E[u](t)=\frac{1}{2\pi}\int_{0}^{2\pi}e[u](t,x)\,dx=E[u](0) & \mbox{(energy)},\\
P[u](t)=\frac{1}{2\pi}\int_0^{2\pi}p[u](t,x)\,dx=P[u](0) & \mbox{(momentum)},
\end{array}
\end{eqnarray}
where  
%\begin{eqnarray*}
%M[u](t)=\int_0^{2\pi}m(u(t,x))\,dx,\quad 
%E[u](t)=\int_{0}^{2\pi}e(u(t,x))\,dx,
%\end{eqnarray*}
%\begin{eqnarray*}
%P[u](t)=\int_0^{2\pi}p(u(t,x))\,dx,
%\end{eqnarray*}
\begin{eqnarray*}%\label{mass-density}
m[u]=|u|^2\quad \mbox{(mass density)},
\end{eqnarray*}
\begin{eqnarray*}%\label{energy-density}
e[u]=\frac{1}{2}|\partial_xu|^2+\frac{\lambda}{4}|u|^2\mathrm{Im}(\overline{u}\partial_xu)+\frac{\lambda^2+2\mu}{12}|u|^6\quad\mbox{(energy density)},
\end{eqnarray*}
\begin{eqnarray*}%\label{momentum-density}
p[u]=-\frac{1}{2}\mathrm{Im}(\overline{u}\partial_xu)-\frac{\lambda}{4}|u|^4\quad \mbox{(momentum density)}.
\end{eqnarray*}
Note that in the case $\mu>-5\lambda^2/16$, the energy $E[u](t)$ assigns positive sign.
In practice if $\mu>-\lambda^2/16$, the absorption of the second term $\frac{\lambda}{4}|u|^2\mathrm{Im}(\overline{u}\partial_xu)$ into the first and third terms yields 
\begin{eqnarray*}
e[u]\ge c\left(|\partial_xu|^2+|u|^6\right),
\end{eqnarray*} 
for some $c>0$, and we conclude
\begin{eqnarray}\label{energy-pos}
E[u](t)\ge c\|\partial_xu(t)\|_{L^2}^2.
\end{eqnarray}
We refer to the cases $\mu>-5\lambda^2/16$ and $\mu<-5\lambda^2/16$ as that the nonlinearity are respectively defocusing and focusing. 
By (\ref{energy-pos}) and $L^2$-norm conservation (\ref{conservation}), in the case $\mu>-5\lambda^2/16$ we have the global well-posedness for any data in $H^s$ for any $s\ge 1$, whereas in the case $\mu\le -5\lambda^2/16$ we require small in $L^2$.

In the present paper, we concentrate on defocusing case as prototype problem.
The purpose of this paper is to investigation the dynamics on resonant clusters for (\ref{dnls}) in the underlying nonlinear energy transfers. 
Notice that such energy exchange processes have been seen with the defocusing cubic nonlinear Schr\"odinger equation on two spacial dimension by Colliander, Keel, Staffilani, Takaoka and Tao \cite{ckstt3}, and with the quintic nonlinear Schr\"odinger equation on one spacial dimension by Gr\'ebert and Thomann \cite{gt1}.
We summarize these results in Remark \ref{2d-cubic} and Remark \ref{1d-quintic}, respectively.

Now it is similar in \cite{gt1}, we begin with a resonant set associated to (\ref{dnls}).
We are interested in whether there exist solutions that have frequencies of the energy exchange observed between the Fourier modes when only four are initially excited. 
\begin{definition}\label{def-res}
%In the sequel, we will use the notation
For $(M,N)\in\mathbb{Z}^2$ with $\lambda(M+N)>0$, we will use the notation
$$
\alpha_1=M,\quad \alpha_2=-3M-N,\quad \beta_1=N,\quad \beta_2=-3N-M,
$$
and $M_*=\max\{|M|,\,|N|\}$.
%From the following discussion, 
With numbers $\alpha_j,~\beta_j~(j=1,2)$ satisfying
\begin{eqnarray}\label{res-res}
\gamma_1+\gamma_3\ne \gamma_2+\gamma_4,
\end{eqnarray}
for $\{\gamma_1,\,\gamma_3\}\ne \{\gamma_2,\,\gamma_4\}$, where each $\gamma_j$ is one of $\alpha_1,\,\alpha_2,\,\beta_1,\,\beta_2$, a resonant set $\Lambda=\Lambda(M,N)$ takes the form
$$
\Lambda=\left\{\alpha_1,\,\alpha_2,\,\beta_1,\,\beta_2\right\}.
$$

Let $P_L$ and $P_H$ be frequency cut-off operators at low and high frequency modes, respectively, defined as follows:
\begin{eqnarray*}
P_Lf(\xi)=\left\{
\begin{array}{ll}
f(\xi), & \mbox{if $|\xi|\lesssim |M+N|$},\\
0, & \mbox{if $|\xi|\gg |M+N|$},
\end{array}
\right.
\end{eqnarray*}
and $P_Hf=(I-P_L)f$.
\end{definition}
\begin{remark}
A typical sufficient condition for (\ref{res-res})  is that $M+N=\mathrm{sgn}\lambda$ and $|M|\gg 1$.
\end{remark}

The main result of this paper is the following theorem.
\begin{theorem}[Case $\lambda=20(M+N)$]\label{main-theorem}
Let $\lambda=20(M+N),~-5\lambda^2/16<\mu\ne 0,~M_*\gg|\lambda|$ and $|\mu|\ll M_*^2$. 
There exist $c_*\in(0,1/2)$ and $T>0$ satisfying $T\to 0$ as $\frac{1}{|\mu|}\to 0$ such that the following holds: there exist $2T$-periodic function $K(t):\mathbb{R}\to (0,1)$ and a solution to (\ref{dnls}) such that $K(0)\le 1/2-c_*$ and $K(T)\ge 1/2+c_*$ and for all $t\in\mathbb{R}$
\begin{eqnarray}\label{formula-dnls}
u(t,x)=\sum_{\xi\in\Lambda}a_{\xi}(t)e^{ix\xi}+e(t,x),\quad e|_{t=0}=0,
\end{eqnarray}
where
\begin{eqnarray*}
|a_{\alpha_{1}}(t)|^2=2|a_{\alpha_{2}}(t)|=K(t),\quad
|a_{\beta_{1}}(t)|^2=2|a_{\beta_{2}}(t)|^2=1-K(t).
\end{eqnarray*}
Moreover for $\delta\in[1/2,1)$ and $|t|\ll \min\{\frac{1}{|\lambda|M_*},\frac{1}{\lambda^2+|\mu|}\}$,
\begin{eqnarray*}
\sum_{\xi\not\in \Lambda}\left|\left(M_*^{\delta}P_L+|\xi|^{\delta}P_H\right){\cal F}_xe(t,\xi)\right|\lesssim \frac{1}{M_*^{\delta}}.
\end{eqnarray*}
\end{theorem}

\begin{remark}
One notices that the periodicity rate $T$ obtained in Theorem \ref{main-theorem} has the bound $\frac{1}{c|\mu|}<T<\frac{c}{|\mu|}$, where the constant $c>1$ is independent of $\lambda,\,\mu,\,M$ and $N$.
Notice that $\frac{1}{|\mu|}\gg\min\{\frac{1}{|\lambda|M_*},\frac{1}{\lambda^2+|\mu|}\}$, provided $|\lambda|M_*\gg |\mu|$ or $|\lambda|\gg 1$.
So it is unclear what the time trajectory of $e(t,x)$ in (\ref{formula-dnls}) describes beyond the time interval $|t|\ll\min\{\frac{1}{|\lambda|M_*},\frac{1}{\lambda^2+|\mu|}\}$ in that setting. 
%If $\lambda^2M_*<|\mu||M+N|$, then $\frac{1}{|\mu|}(\frac{\lambda}{M+N})^2<\frac{1}{|M+N|M_*}$ holds.
\end{remark}

\begin{remark}
The above theorem gives the approximation
$$
u(t,x)\approx \sum_{\xi\in\Lambda}a_{\xi}(t)e^{ix\xi},
$$
for $|t|\ll \min\{\frac{1}{|\lambda|M_*},\frac{1}{\lambda^2+|\mu|}\}$, where masses $a_{\beta_j}$ are excited initially compared these of $a_{\alpha_j}$, then move decreasingly to that of $a_{\alpha_j}$, in periodically.
%So one can view this example as illustrating the exchange of energy  
\end{remark}

\begin{remark}\label{2d-cubic}
For the defocusing cubic nonlinear Schr\"odinger equation on two spatial dimension:
$$
i\partial_tu+\Delta u=|u|^2u,
$$
one can obtain the low-to-high energy cascade scenario, see \cite{ckstt3}.
In the proof, an instability phenomenon of Arnold diffusion phenomena is considered in perturbed Hamiltonian systems.    
Indeed, for any $s>1$ and $K\gg 1\gg\varepsilon$ there exist a global solution $u(t)$ of 2d-NLS on $\mathbb{T}^2$ and time $T$ satisfying that
$$
\|u(0)\|_{H^s}\le \varepsilon,\quad \|u(T)\|_{H^s}\ge K.
$$
\end{remark}

\begin{remark}\label{1d-quintic}
In \cite{gt1}, one striking example of a similar situation is considered for the quintic nonlinear Schr\"odinger equation on one spacial dimension:
\begin{eqnarray}\label{quintic-nls}
i\partial_tu+\partial_x^2u=\nu|u|^4u.
\end{eqnarray}
For $(k,\,n)\in\mathbb{Z}^2$ with $k\ne 0$, a resonant set $\Lambda_1$ corresponding to (\ref{quintic-nls}) is the following:
$$
\Lambda_1=\left\{\,\alpha'_1,\,\alpha_2',\,\beta'_1,\,\beta'_2\,\right\},
$$
where
$$
(\alpha'_1,\,\alpha_1',\,\beta_1',\,\beta_2')=(n+3k,\,n,\,n+k,\,n+4k).
$$
Comparing to the result stated in Theorem \ref{main-theorem}, they obtained dynamical properties for (\ref{quintic-nls}) in stronger sense.
More precisely, it is obtained that there exist $T>0,~\nu_0>0,~2T$-periodic function $K:\mathbb{R}\to (0,1)$ with $K(0)\le 1/4$ and $K(T)\ge 3/4$ so that for all $0<\nu<\nu_0$ and $s\in \mathbb{R}$, there exists a solution of (\ref{quintic-nls})  satisfying
\begin{eqnarray}\label{formula-nls}
u(t)=\sum_{\xi\in \Lambda_1}a_{\xi}(t)e^{ix\xi}+\left(\nu^{1/4}+t\nu^{3/2}\right)O_{H^s}(1)
\end{eqnarray}
for $|t|\lesssim \nu^{-3/2}$, where
$$
|a_{n+3k}(t)|^2=2|a_{n}(t)|=K(\nu t),
$$
$$
|a_{n+k}(t)|^2=2|a_{n+4k}(t)|^2=1-K(\nu t).
$$
The proof is centered on the resonant normal form approach.
The formula (\ref{formula-dnls}) looks very similar to that of (\ref{formula-nls}), but this conveys different information because of obstructions by the derivative nonlinearity in (\ref{dnls}).
%is a slightly different. route , which is essential part of the nonlinear  obstruction
\end{remark}

%\begin{remark}\label{scaling}
We now turn to equation adjusting the amplitude scaling:
$$
u(t,x)\mapsto v(t,x)=\sqrt{\frac{20(M+N)}{\lambda'}}u(t,x),\quad\lambda'(M+N)>0,
$$
which transforms the equation (\ref{dnls}) in the case $\lambda=20(M+N)$ to the following:
\begin{eqnarray}\label{sca-dnls}
i\partial_tv+\partial_x^2v=-\lambda'v^2\partial_x\overline{v}+\mu'|v|^4v,
\end{eqnarray}
where
$$
\mu'=\mu\left(\frac{20(M+N)}{\lambda'}\right)^2.
$$
%\end{remark}
Apply this to Theorem \ref{main-theorem}, we obtain the immediate corollary, which is equivalent to the counterpart of Theorem \ref{main-theorem}.
%In the subsequent sections, we show Corollary \ref{main-corollary} instead of Theorem \ref{main-theorem}.
\begin{corollary}[Case $\lambda'(M+N)>0$]\label{main-corollary}
Let $-5\lambda'^2/16<\mu'\ne 0,~M_*\gg |M+N|$ and $ M_*|M+N|\gg |\lambda'|\sqrt{|\mu'|}$.
There exist $c_*\in(0,1/2)$ and $T>0$ satisfying $T\to 0$ as $\frac{1}{|\mu'|}(\frac{M+N}{\lambda'})^2\to 0$ such that the following holds: there exist $2T$-periodic function $K(t):\mathbb{R}\to (0,1)$ and a global solution to (\ref{sca-dnls}) such that $K(0)\le 1/2-c_*$ and $K(T)\ge 1/2+c_*$ and for all $t\in\mathbb{R}$
\begin{eqnarray*}%\label{formula-dnls}
v(t,x)=\sqrt{\frac{20(M+N)}{\lambda'}}\left(\sum_{\xi\in\Lambda}a_{\xi}(t)e^{ix\xi}+e(t,x)\right),\quad e|_{t=0}=0,
\end{eqnarray*}
where
\begin{eqnarray*}
|a_{\alpha_{1}}(t)|^2=2|a_{\alpha_{2}}(t)|=K(t),\quad
|a_{\beta_{1}}(t)|^2=2|a_{\beta_{2}}(t)|^2=1-K(t).
\end{eqnarray*}
Moreover for $\delta\in[1/2,1)$ and $|t|\ll \min\{\frac{1}{|M+N|M_*},\frac{(M+M)^2}{(M+N)^4+\lambda'^2|\mu'|}\}$,
\begin{eqnarray*}
\sum_{\xi\not\in \Lambda}\left|\left(M_*^{\delta}P_L+|\xi|^{\delta}P_H\right){\cal F}_xe(t,\xi)\right|\lesssim \frac{1}{M_*^{\delta}}.
\end{eqnarray*}
\end{corollary}

Theorem \ref{main-theorem} is obtained by analyzing an infinitely dimensional system of ordinary differential equations that include the resonant components on nonlinearity of (\ref{dnls}).
We first make the Fourier transform to the equation (\ref{dnls}) and build a resonant approximation model (toy model), which is an ordinary differential system.
This is accomplished by investigating the long time dynamics of toy model.
The error estimate is proved to approximate well the solution of (\ref{dnls}) by that of the toy model by using the bootstrap argument.

\begin{notation}
We prefer to use the notation $\langle\cdot\rangle=(1+|\cdot|^2)^{1/2}$.
The over dot $\dot{a}(t)$ denotes the time derivative of $a(t)$. 

The Fourier transform with respect to the space variable (discrete Fourier transform) is defined by
$$
{\cal F}_xf(\xi)=\frac{1}{\sqrt{2\pi}}\int_0^{2\pi}e^{-ix\xi}f(x)\,dx,\quad \xi\in\mathbb{Z}.
$$
We abbreviate ${\cal F}_xf$ as $\widehat{f}$. 

For $1\le p,q\le \infty$, we use the mixed norm notation $\|f\|_{L^q_tL^p_x}$ with norm
$$
\|f\|_{L^q_tL^p_x}=\left(\int_{-\infty}^{\infty}\|f(t)\|_{L^p_x(\mathbb{T})}^q\,dt\right)^{1/q},
$$
with the obvious modification when $q=\infty$.
For $T>0$, we also use $\|f\|_{L^q_TL^p_x}$ to denote the norm
$$
\|f\|_{L^q_TL^p_x}=\left(\int_0^T\|f(t)\|_{L^p_x(\mathbb{T})}^q\,dt\right)^{1/q},
$$
with the obvious modification when $q=\infty$.
We also use $L_t^p\ell_{\xi}^q$ in a similar manner.

We use $c,~C$ to denote various constants that does not depend on $\lambda,\,\mu,\,M$ and $N$.
We use $A\lesssim B$ to denote $A\le CB$ for some constant $C>0$.
Similarly, we write $A\ll B$ to mean $A\le c B$ for some small constant $c>0$.

For complex-valued $n$ functions $f_1,~f_2,\ldots,~f_n$ defined on the set $\mathbb{Z}$ of integers, we write
the discrete convolution (convolution sum) $[f_1*f_2*\ldots *f_n](\xi)$ as
$$
[f_1*f_2*\ldots *f_n](\xi)=\sum_{*}\prod_{j=1}^nf_1(\xi_j),
$$
where $\sum_*$ denotes a summation over the set where $\xi_1+\xi_2+\ldots +\xi_n=\xi$.
%Also write
%$$
%[g_1*g_2*\ldots*g_n](\tau)=\int_*\prod_{j=1}^ng_j(\tau_j),
%$$
%where $\int_*$ denotes an integration over the set where $\tau_1+\tau_2+\ldots+\tau_n=\tau$.
\end{notation}

The rest of this paper is organized as follows.
In Section \ref{fourieransatz}, we reduce the equation (\ref{dnls}) to that for the Fourier coefficients, by taking the advantage of conservation laws (\ref{conservation}). 
In Section \ref{toymodel}, we construct an appropriate ordinary differential equation (ODE) system as theoretical toy model for (\ref{dnls}).
In Section \ref{odedynamics}, we solve the ODE system.
In Section \ref{error-sec}, the equation that we derived in the previous section should approximate the dynamics of the first differential equation (\ref{dnls}).

\section{Fourier Ansatz}\label{fourieransatz}
\indent

In this section we derive the resonant approximation model related to (\ref{dnls}).
We consider the smooth solutions thought the paper, and write
\begin{eqnarray}\label{four-ansatz}
u(t,x)=\sum_{\xi\in\mathbb{Z}}a_{\xi}(t)e^{i(x\xi+t\xi^2)}.
\end{eqnarray}
Substituting this into the equation (\ref{dnls}) and using
\begin{eqnarray}\label{identity}
\xi_1^2-\xi_2^2+\xi_3^2-\xi^2=2(\xi_1-\xi_2)(\xi_1-\xi)=2(\xi_3-\xi_2)(\xi_3-\xi),
\end{eqnarray}
for $\xi=\xi_1-\xi_2+\xi_3$, we obtain the equation for $a_{\xi}=a_{\xi}(t)$:
\begin{eqnarray}\label{f-1-dnls}
i\dot{a}_{\xi}& = & \lambda\xi|a_{\xi}|^2a_{\xi}\nonumber\\
& & -\lambda\sum_{\scriptstyle * \atop{\scriptstyle \{\xi_1,\xi_3\}\cap\{\xi_2,\xi\}=\emptyset}}\xi_2a_{\xi_1}\overline{a_{\xi_2}}a_{\xi_3}e^{2it(\xi_1-\xi_2)(\xi_1-\xi)}\nonumber\\
& & -2\lambda\left(\sum_{\xi'\in\mathbb{Z}}\xi'|a_{\xi'}|^2\right)a_{\xi}\nonumber\\
& & +\mu\sum_* a_{\xi_1}\overline{a_{\xi_2}}a_{\xi_3}\overline{a_{\xi_4}}a_{\xi_5}e^{it(\xi_1^2-\xi_2^2+\xi_3^2-\xi_4^2+\xi_5^2-\xi^2)},
\end{eqnarray}
by splitting the cubic convolution sum on $\xi=\xi_1-\xi_2+\xi_3$ into two types of sums as same as in \cite{ta2}; $\{\xi_1,\xi_3\}\cap \{\xi_2,\xi\}=\emptyset$ and $\{\xi_1,\xi_3\}\cap \{\xi_2,\xi\}\ne\emptyset$.
%$$
%\sum_{*}=-\sum_{\scriptstyle * \atop{\scriptstyle \xi_1=\xi_2=\xi_3}}+\sum_{\scriptstyle * \atop{\scriptstyle \{\xi_1,\xi_3\}\cap\{\xi_2,\xi\}=\emptyset}}+\sum_{\scriptstyle * \atop{\scriptstyle \xi_1=\xi_2}}+\sum_{\scriptstyle * \atop{\scriptstyle \xi_2=\xi_3}}.
%$$

Denote $P_0=P[u](0)$ and $M_0=M[u](0)$.
Taking the momentum conservation (\ref{conservation}) via Fourier transformation, we obtain
$$
\sum_{\xi'\in\mathbb{Z}}\xi'|a_{\xi'}|^2=-\frac{\lambda}{2}\sum_{\xi_1+\xi_3=\xi_2+\xi_4}a_{\xi_1}\overline{a_{\xi_2}}a_{\xi_3}\overline{a_{\xi_4}}e^{2it(\xi_1-\xi_2)(\xi_1-\xi_4)}-2P_0.
$$
On substituting this into the second term on the right-hand side of (\ref{f-1-dnls}), we have the following form
\begin{eqnarray}\label{f-2-dnls}
i\dot{a}_{\xi} & = & \lambda\xi|a_{\xi}|^2a_{\xi}+4\lambda P_0a_{\xi}\nonumber\\
& & -\lambda\sum_{\scriptstyle * \atop{\scriptstyle \{\xi_1,\xi_3\}\cap\{\xi_2,\xi\}=\emptyset}}\xi_2a_{\xi_1}\overline{a_{\xi_2}}a_{\xi_3}e^{2it(\xi_1-\xi_2)(\xi_1-\xi)}\nonumber\\
& & +\lambda^2a_{\xi}\sum_{\xi_1+\xi_3=\xi_2+\xi_4}a_{\xi_1}\overline{a_{\xi_2}}a_{\xi_3}\overline{a_{\xi_4}}e^{2it(\xi_1-\xi_2)(\xi_1-\xi_4)}\nonumber\\
& & +\mu\sum_* a_{\xi_1}\overline{a_{\xi_2}}a_{\xi_3}\overline{a_{\xi_4}}a_{\xi_5}e^{it(\xi_1^2-\xi_2^2+\xi_3^2-\xi_4^2+\xi_5^2-\xi^2)}.
\end{eqnarray}
Note that $(\xi_1-\xi_2)(\xi_1-\xi)\ne 0$ when $\{\xi_1,\xi_3\}\cap\{\xi_2,\xi\}=\emptyset$ in the third term on the right-hand side of (\ref{f-2-dnls}).
When $(\xi_1-\xi_2)(\xi_1-\xi_4)\ne 0$, it oscillates in space-time at different frequencies, which illustrates hidden smoothing effect.
Thus for the third term on the right-hand side of (\ref{f-2-dnls}), we divide the sum into two cases: $\{\xi_1,\xi_3\}\cap\{\xi_2,\xi\}=\emptyset$ and $\{\xi_1,\xi_3\}\cap\{\xi_2,\xi\}\ne\emptyset$, then
\begin{eqnarray*}
& & \sum_{\xi_1+\xi_3=\xi_2+\xi_4}a_{\xi_1}\overline{a_{\xi_2}}a_{\xi_3}\overline{a_{\xi_4}}e^{2it(\xi_1-\xi_2)(\xi_1-\xi_4)} \\
&  = & 2\left(\sum_{\xi'\in\mathbb{Z}}|a_{\xi'}|^2\right)^2-\sum_{\xi'\in\mathbb{Z}}|a_{\xi'}|^4\\
& & +\sum_{\scriptstyle \xi_1+\xi_3=\xi_2+\xi_4 \atop{\scriptstyle \{\xi_1,\xi_3\}\cap\{\xi_2,\xi\}=\emptyset} }a_{\xi_1}\overline{a_{\xi_2}}a_{\xi_3}\overline{a_{\xi_4}}e^{2it(\xi_1-\xi_2)(\xi_1-\xi_4)}\\
& = &  2M_0^2-\sum_{\xi'\in\mathbb{Z}}|a_{\xi'}|^4 +\sum_{\scriptstyle \xi_1+\xi_3=\xi_2+\xi_4 \atop{\scriptstyle \{\xi_1,\xi_3\}\cap\{\xi_2,\xi\}=\emptyset} }a_{\xi_1}\overline{a_{\xi_2}}a_{\xi_3}\overline{a_{\xi_4}}e^{2it(\xi_1-\xi_2)(\xi_1-\xi_4)}.
\end{eqnarray*}
We also divide the sum in the last term on the right-hand side of (\ref{f-2-dnls}) into two cases; $\{\xi_1,\,\xi_3,\,\xi_5\}=\{\xi_2,\,\xi_4,\,\xi\}$ and $\{\xi_1,\,\xi_3,\,\xi_5\}\ne \{\xi_2,\,\xi_4,\,\xi\}$.
On the set $\{\xi_1,\,\xi_3,\,\xi_5\}=\{\xi_2,\,\xi_4,\,\xi\}$, we split this case into three disjoint subsets; (i) $\xi_1=\xi,\,\{\xi_3,\,\xi_5\}=\{\xi_2,\,\xi_4\}$, (ii) $\xi_1\ne \xi,\,\{\xi_1,\,\xi_5\}=\{\xi_2,\,\xi_4\}$, (iii) $\xi_1\ne \xi,\,\xi_3\ne \xi,\,\{\xi_1,\,\xi_3\}=\{\xi_2,\,\xi_3\}$.
The contribution of the subset (i) to the last term on the right-hand side of (\ref{f-2-dnls}) is equal to
\begin{eqnarray*}
& & a_{\xi}\left(2\sum_{\xi'\in\mathbb{Z}}|a_{\xi'}|^2\sum_{\xi''\ne \xi'}|a_{\xi''}|^2+\sum_{\xi'\in\mathbb{Z}}|a_{\xi'}|^4\right)\\
& = & a_{\xi}\left(2\left(\sum_{\xi'\in\mathbb{Z}}|a_{\xi'}|^2\right)^2-\sum_{\xi'\in\mathbb{Z}}|a_{\xi'}|^4\right)\\
& = & a_{\xi}\left(2M_0^2-\sum_{\xi'\in\mathbb{Z}}|a_{\xi'}|^4\right).
\end{eqnarray*}
The contribution of the subset (ii) to the last term on the right-hand side of (\ref{f-2-dnls}) is equal to
\begin{eqnarray*}
& & a_{\xi}\left(2\sum_{\xi'\ne \xi}|a_{\xi'}|^2\sum_{\xi''\ne \xi'}|a_{\xi''}|^2+\sum_{\xi'\ne \xi}|a_{\xi'}|^2\right)\\
& = & a_{\xi}\left(\left(\sum_{\xi'\in\mathbb{Z}}|a_{\xi'}|^2\right)^2-2|a_{\xi}|^2\sum_{\xi'\in\mathbb{Z}}|a_{\xi'}|^2-\sum_{\xi'\in\mathbb{Z}}|a_{\xi'}|^4+|a_{\xi}|^4\right)\\
& = & a_{\xi}\left(2M_0^2-2|a_{\xi}|^2M_0-\sum_{\xi'\in\mathbb{Z}}|a_{\xi'}|^4+|a_{\xi}|^4\right)
.
\end{eqnarray*}
The contribution of the subset (iii) to the last term on the right-hand side of (\ref{f-2-dnls}) is equal to
\begin{eqnarray*}
& & a_{\xi}\left(2\sum_{\xi'\ne \xi}|a_{\xi'}|^2\sum_{\xi''\not\in\{\xi,\,\xi'\}}|a_{\xi''}|^2+\sum_{\xi'\ne \xi}|a_{\xi'}|^2\right)\\
& = & a_{\xi}\left(2\left(\sum_{\xi'\in\mathbb{Z}}|a_{\xi'}|^2\right)^2-4|a_{\xi}|^2\sum_{\xi'\in\mathbb{Z}}|a_{\xi'}|^2-\sum_{\xi'\in\mathbb{Z}}|a_{\xi'}|^4+3|a_{\xi}|^4\right)\\
& = & a_{\xi}\left(2M_0^2-4|a_{\xi}|^2M_0-\sum_{\xi'\in\mathbb{Z}}|a_{\xi'}|^4+3|a_{\xi}|^4\right).
\end{eqnarray*}
Collecting the above results, we have that
\begin{eqnarray*}
& & \sum_* a_{\xi_1}\overline{a_{\xi_2}}a_{\xi_3}\overline{a_{\xi_4}}a_{\xi_5}e^{it(\xi_1^2-\xi_2^2+\xi_3^2-\xi_4^2+\xi_5^2-\xi^2)}\\
& = & \left(6M_0^2-6M_{0}|a_{\xi}|^2-3\sum_{\xi'\in\mathbb{Z}}|a_{\xi'}|^4+4|a_{\xi}|^4\right)a_{\xi}\\
& & +\sum_{\scriptstyle * \atop{\scriptstyle \{\xi_1,\,\xi_3,\,\xi_5\}\ne \{\xi_2,\,\xi_4,\,\xi\}}}a_{\xi_1}\overline{a_{\xi_2}}a_{\xi_3}\overline{a_{\xi_4}}a_{\xi_5}e^{it(\xi_1^2-\xi_2^2+\xi_3^2-\xi_4^2+\xi_5^2-\xi^2)}.
\end{eqnarray*}
From this, the formula (\ref{f-2-dnls}) is described as
\begin{eqnarray}\label{f-dnls}
i\dot{a}_{\xi} & = & (\lambda\xi-6M_0\mu)|a_{\xi}|^2a_{\xi}+2\left(M_0^2(\lambda^2+3\mu)+2\lambda P_0\right)a_{\xi}\nonumber\\
& & +4\mu|a_{\xi}|^4a_{\xi}-(\lambda^2+3\mu)\left(\sum_{\xi'\in\mathbb{Z}}|a_{\xi'}|^4\right)a_{\xi}\nonumber\\
& & -\lambda\sum_{\scriptstyle * \atop{\scriptstyle \{\xi_1,\xi_3\}\cap\{\xi_2,\xi\}=\emptyset}}\xi_2a_{\xi_1}\overline{a_{\xi_2}}a_{\xi_3}e^{2it(\xi_1-\xi_2)(\xi_1-\xi)}\nonumber\\
& & +\lambda^2a_{\xi}\sum_{\scriptstyle \xi_1+\xi_3=\xi_2+\xi_4 \atop{\scriptstyle \{\xi_1,\xi_3\}\cap\{\xi_2,\xi\}=\emptyset} }a_{\xi_1}\overline{a_{\xi_2}}a_{\xi_3}\overline{a_{\xi_4}}e^{2it(\xi_1-\xi_2)(\xi_1-\xi_4)}\nonumber\\
& & +\mu\sum_{\scriptstyle * \atop{\scriptstyle \{\xi_1,\,\xi_3,\,\xi_5\}\ne \{\xi_2,\,\xi_4,\,\xi\}}} a_{\xi_1}\overline{a_{\xi_2}}a_{\xi_3}\overline{a_{\xi_4}}a_{\xi_5}e^{it(\xi_1^2-\xi_2^2+\xi_3^2-\xi_4^2+\xi_5^2-\xi^2)}.
\end{eqnarray} 
The cubic interactions through resonance in the first term on the right-hand side of (\ref{f-dnls})  are self-interaction, which enhance the phase oscillations.
%We now consider the small cardinal 
In next section, we seek the contribution of quintic resonant interactions to the last term on the right-hand side of (\ref{f-dnls}).

\section{Toy model}\label{toymodel}
\indent

In this section we shall construct resonant, finite dimensional truncation model to (\ref{f-dnls}).
%By Remark \ref{scaling}, 
Let
$$
\lambda=20(\alpha_1+\beta_1)=20(M+N).
$$
We start out by referring to the observation on the resonant set used in \cite{gt1}, which can be helpful in predicting the resonant set related to the equation (\ref{f-dnls}).
\begin{lemma}[{Gr\'ebert-Thomann \cite[Lemmas 2.1 and 2.2]{gt1}}]
Assume that $\xi_j\in\mathbb{Z}~(1\le j\le 6)$ satisfy
$$
\xi_1+\xi_3+\xi_5=\xi_2+\xi_4+\xi_6,
$$
$$
\xi_1^2+\xi_3^2+\xi_5^2=\xi_2^2+\xi_4^2+\xi_6^2,
$$
and $\{\xi_1,\,\xi_3,\,\xi_5\}\ne \{\xi_2,\,\xi_4,\,\xi_6\}$.
Then $\{\xi_1,\,\xi_3,\,\xi_5\}\cap\{\xi_2,\,\xi_4,\,\xi_6\}=\emptyset$.
Moreover the cardinal of $\{\xi_j\mid 1\le j\le 6\}$ is greater than or equal to $4$.
\end{lemma}

Considering the influence by  the first term on the right-hand side of (\ref{f-dnls}), we define the resonant interaction set
$$
\Lambda^*=\left\{\{\alpha_1,\,\alpha_1,\,\alpha_2\},\,\{\beta_1,\,\beta_1,\,\beta_2\}\right\},
$$
where $\alpha_j$ and $\beta_j$ are defined in Definition \ref{def-res}.
We will use the notation
$$
\left\{\{\xi_1,\,\xi_3,\,\xi_5\},\,\{\xi_2,\,\xi_4,\,\xi\}\right\}=\Lambda^*
$$
to represent that two of $\xi_1,\,\xi_3,\,\xi_5$ are $\alpha_1$, one of $\xi_1,\,\xi_3,\,\xi_5$ are $\alpha_2$, two of $\xi_2,\,\xi_4,\,\xi$ are $\beta_1$, one of $\xi_2,\,\xi_4,\,\xi$ is $\beta_2$, and its symmetrical objects.

Having discussed in Section \ref{fourieransatz}, we define the resonant truncation ODE model of (\ref{f-dnls}) by
\begin{eqnarray}\label{m-dnls}
i\dot{c}_{\xi} & = & (\lambda\xi-6M_0\mu)|c_{\xi}|^2c_{\xi}+2\left(M_0^2(\lambda^2+3\mu)+2\lambda P_0\right)c_{\xi}\nonumber\\
& & +4\mu|c_{\xi}|^4c_{\xi}-(\lambda^2+3\mu)\left(\sum_{\xi'\in\Lambda}|c_{\xi'}|^4\right)c_{\xi}\nonumber\\
& & +\mu\sum_{\left\{\{\xi_1,\,\xi_3,\,\xi_5\},\,\{\xi_2,\,\xi_4,\,\xi\}\right\}=\Lambda^*} c_{\xi_1}\overline{c_{\xi_2}}c_{\xi_3}\overline{c_{\xi_4}}c_{\xi_5}e^{it(\xi_1^2-\xi_2^2+\xi_3^2-\xi_4^2+\xi_5^2-\xi^2)},
\end{eqnarray}
with initial data
\begin{eqnarray}\label{md-dnls}
c_{\xi}(0)=a_{\xi}(0),
\end{eqnarray}
where $\xi\in \Lambda$.

One of the most important features of the above ODE system is the following lemma.
\begin{lemma}\label{integral}
Let $\{c_{\xi}\}_{\xi\in \Lambda}$ be a solution to (\ref{m-dnls})-(\ref{md-dnls}) on the time interval $I$.
Then we have the relations:
\begin{eqnarray}\label{alpha1-beta1}
\frac{d}{dt}\left(|c_{\alpha_1}(t)|^2+|c_{\beta_1}(t)|^2\right)=0,
\end{eqnarray}
\begin{eqnarray}\label{alpha1-beta2}
\frac{d}{dt}\left(|c_{\alpha_2}(t)|^2+|c_{\beta_2}(t)|^2\right)=0,
\end{eqnarray}
\begin{eqnarray}\label{alpha1-alpha2}
\frac{d}{dt}\left(|c_{\alpha_1}(t)|^2-2|c_{\alpha_2}(t)|^2\right)=0,
\end{eqnarray}
\begin{eqnarray}\label{beta1-beta2}
\frac{d}{dt}\left(|c_{\beta_1}(t)|^2-2|c_{\beta_2}(t)|^2\right)=0,
\end{eqnarray}
for $t\in I$.
\end{lemma}
%\begin{remark} 
%
%
%By the Arnold-Liouville theorem, the above lemma in fact ensures that there is four independent asserts that the system (\ref{m-dnls})-(\ref{md-dnls}) is completely integrable.
%
%\end{remark}
\noindent
{\it Proof.}
The proof is straightforward calculation.
We only give the proof of (\ref{alpha1-beta1}) and (\ref{alpha1-alpha2}), as the proofs for (\ref{alpha1-beta2}) and (\ref{beta1-beta2}) are analogous.

For (\ref{alpha1-beta1}) we multiply equation (\ref{m-dnls})  for $\xi=\alpha_1$ by $2\overline{c_{\alpha_1}}$ and take the imaginary part to get
\begin{eqnarray}\label{c_1}
\frac{d}{dt}|c_{\alpha_1}|^2& = & 2\mu\Im\sum_{\scriptstyle \{\xi_1,\,\xi_3,\,\xi_5\}=\{\beta_1,\,\beta_1,\,\beta_2\}\atop{\scriptstyle \{\xi_2,\,\xi_4\}=\{\alpha_1,\,\alpha_2\}}} c_{\xi_1}\overline{c_{\xi_2}}c_{\xi_3}\overline{c_{\xi_4}}c_{\xi_5}\overline{c_{\alpha_1}}e^{it(\xi_1^2-\xi_2^2+\xi_3^2-\xi_4^2+\xi_5^2-\alpha_1^2)}\nonumber\\
& = & 2^3\mu \Im c_{\beta_1}\overline{c_{\alpha_2}}c_{\beta_3}\overline{c_{\alpha_2}}c_{\beta_2}\overline{c_{\alpha_1}}e^{it(2\beta_1^2+\beta_2^2-2\alpha_1^2-\alpha_2^2)},
\end{eqnarray}
where we use that
$$
\xi_1^2-\xi_2^2+\xi_3^2-\xi_4^2+\xi_5^2-\alpha_1^2=2\beta_1^2+\beta_2^2-2\alpha_1^2-\alpha_2^2,
$$
for $\{\xi_1,\,\xi_3,\,\xi_5\}=\{\beta_1,\,\beta_1,\,\beta_2\}$ and $\{\xi_2,\,\xi_4\}=\{\alpha_1,\,\alpha_2\}$.
A similar argument as above but using equation (\ref{m-dnls}) for $\xi=\beta_1$ yields
\begin{eqnarray}\label{c_2}
\frac{d}{dt}|c_{\beta_1}|^2 & = & 2\mu\Im\sum_{\scriptstyle \{\xi_1,\,\xi_3,\,\xi_5\}=\{\alpha_1,\,\alpha_1,\,\alpha_2\}\atop{\scriptstyle \{\xi_2,\,\xi_4\}=\{\beta_1,\,\beta_2\}}} c_{\xi_1}\overline{c_{\xi_2}}c_{\xi_3}\overline{c_{\xi_4}}c_{\xi_5}\overline{c_{\beta_1}}e^{it(\xi_1^2-\xi_2^2+\xi_3^2-\xi_4^2+\xi_5^2-\beta_1^2)}\nonumber\\
 & = & 2^3\mu\Im c_{\alpha_1}\overline{c_{\beta_1}}c_{\alpha_1}\overline{c_{\beta_2}}c_{\alpha_2}\overline{c_{\beta_1}}e^{-it(2\beta_1^2+\beta_2^2-2\alpha_1^2-\alpha_2^2)}.
\end{eqnarray}
Collecting the information in (\ref{c_1}) and (\ref{c_2}) we can obtain
\begin{eqnarray*}
\frac{d}{dt}\left(|c_{\alpha_1}|^2+|c_{\beta_1}|^2\right) & = & 2^3\mu\Im\left( c_{\beta_1}\overline{c_{\alpha_2}}c_{\beta_3}\overline{c_{\alpha_2}}c_{\beta_2}\overline{c_{\alpha_1}}e^{it(2\beta_1^2+\beta_2^2-2\alpha_1^2-\alpha_2^2)}\right.\\
& & +\left.c_{\alpha_1}\overline{c_{\beta_1}}c_{\alpha_1}\overline{c_{\beta_2}}c_{\alpha_2}\overline{c_{\beta_1}}e^{-it(2\beta_1^2+\beta_2^2-2\alpha_1^2-\alpha_2^2)}\right)\\
& = & 0.
\end{eqnarray*}

On the other hand, to prove (\ref{alpha1-alpha2}), we will see that
\begin{eqnarray*}
\frac{d}{dt}|c_{\alpha_2}|^2& = & 2\mu\Im\sum_{\scriptstyle \{\xi_1,\,\xi_3,\,\xi_5\}=\{\beta_1,\,\beta_1,\,\beta_2\}\atop{\scriptstyle \{\xi_2,\,\xi_4\}=\{\alpha_1,\,\alpha_1\}}} c_{\xi_1}\overline{c_{\xi_2}}c_{\xi_3}\overline{c_{\xi_4}}c_{\xi_5}\overline{c_{\alpha_2}}e^{it(\xi_1^2-\xi_2^2+\xi_3^2-\xi_4^2+\xi_5^2-\alpha_2^2)}\nonumber\\
& = & 2^2\mu \Im c_{\beta_1}\overline{c_{\alpha_1}}c_{\beta_3}\overline{c_{\alpha_1}}c_{\beta_2}\overline{c_{\alpha_2}}e^{it(2\beta_1^2+\beta_2^2-2\alpha_1^2-\alpha_2^2)}.
\end{eqnarray*}
Hence the argument used in above gives
$$
\frac{d}{dt}\left(|c_{\alpha_1}|^2-2|c_{\alpha_2}|^2\right)=0,
$$
which concludes the proof of (\ref{alpha1-alpha2}).
\qed

\begin{remark}\label{ode-global}
The conservation laws in (\ref{alpha1-beta1})-(\ref{alpha1-alpha2}) are used to the study of properties of the system of (\ref{m-dnls}) including the global existence results.
In particular, the ODE system (\ref{m-dnls})-(\ref{md-dnls}) has global solutions.
We will show that the system has periodic orbits in times.
\end{remark}

Let us set data $\{a_{\xi}(0)\}$ in (\ref{md-dnls}) where
\begin{eqnarray*}
|a_{\alpha_1}(0)|^2+|a_{\beta_1}(0)|^2=2|a_{\alpha_2}(0)|^2+2|a_{\beta_2}(0)|^2=1,
\end{eqnarray*}
\begin{eqnarray*}
|a_{\alpha_1}(0)|^2-2|a_{\alpha_2}(0)|^2= |a_{\beta_1}(0)|^2-2|a_{\beta_2}(0)|^2=0.
\end{eqnarray*}
By Lemma \ref{integral} and Remark \ref{ode-global}, we have that for all $t\in\mathbb{R}$,
\begin{eqnarray}\label{c_3}
|c_{\alpha_1}(t)|^2+|c_{\beta_1}(t)|^2=2|c_{\alpha_2}(t)|^2+2|c_{\beta_2}(t)|^2=1,
\end{eqnarray}
\begin{eqnarray}\label{c_4}
|c_{\alpha_1}(t)|^2-2|c_{\alpha_2}(t)|^2= |c_{\beta_1}(t)|^2-2|c_{\beta_2}(t)|^2=0.
\end{eqnarray}
We now write
\begin{eqnarray*}
|c_{\alpha_1}(t)|^2=K(t),\quad |c_{\alpha_2}|^2=\frac{K(t)}{2},
\end{eqnarray*}
\begin{eqnarray*}
|c_{\beta_1}(t)|^2=1-K(t),\quad |c_{\beta_2}|^2=\frac{1-K(t)}{2}.
\end{eqnarray*}
\begin{remark}\label{km}
From the relation in (\ref{c_3}) and (\ref{c_4}), we easily see that
\begin{eqnarray}\label{q}
M_0=\frac{3}{2},\quad \sum_{\xi'\in\Lambda}|c_{\xi'}(t)|^4=\frac{5}{4}-\frac{5}{2}K(t)\left(1-K(t)\right).
\end{eqnarray}
\end{remark}
The equation (\ref{m-dnls}) has a gauge freedom.
On writing
$$
d_{\xi}(t)=c_{\xi}(t)\exp \left(iG_{\xi}(t)\right),
$$
where
$$
\dot{G}_{\xi}(t)=\lambda\xi|c_{\xi}(t)|^2+2\left(M_0^2(\lambda^2+3\mu)+2\lambda P_0\right)-(\lambda^2+3\mu)\sum_{\xi'\in\Lambda}|c_{\xi'}(t)|^4,
$$
$$
G_{\xi}(0)=0,
$$
and using (\ref{q}), we have that the equation (\ref{m-dnls}) becomes the following equation for $d_{\xi}$,
\begin{eqnarray}\label{d-dnls}
i\dot{d}_{\xi} & = & \mu |d_{\xi}|^2(4|d_{\xi}|^2-9)d_{\xi}\nonumber\\
& & +\mu\sum_{\left\{\{\xi_1,\,\xi_3,\,\xi_5\},\,\{\xi_2,\,\xi_4,\,\xi\}\right\}=\Lambda^*} d_{\xi_1}\overline{d_{\xi_2}}d_{\xi_3}\overline{d_{\xi_4}}d_{\xi_5}e^{i{\cal I}},
\end{eqnarray}
where ${\cal I}(0)=0$ and 
\begin{eqnarray*}%\label{phase-1}
\dot{\cal I}(t) & = & \xi_1^2-\xi_2^2+\xi_3^2-\xi_4^2+\xi_5^2-\xi^2\nonumber\\
& & -\lambda\left(\xi_1|d_{\xi_1}(t)|^2-\xi_2|d_{\xi_2}(t)|^2+\xi_3|d_{\xi_3}(t)|^2-\xi_4|d_{\xi_4}(t)|^2+\xi_5|d_{\xi_5}(t)|^2-\xi|d_{\xi}(t)|^2\right).
\end{eqnarray*}
Let us discuss the contribution of the right-hand of (\ref{d-dnls}).

\begin{lemma}\label{phase-2}
Let %$\lambda=-20(M+N),~(\alpha_1,\,\alpha_2,\,\beta_1,\,\beta_2)=(M,\,-3M-N,\,N,\,-3N-M)$,
$$
|d_{\alpha_1}(t)|^2=K(t),\quad |d_{\alpha_2}(t)|^2=\frac{K(t)}{2},
$$
$$
|d_{\beta_1}(t)|^2=1-K(t),\quad |d_{\beta_2}(t)|^2=\frac{1-K(t)}{2},
$$
and $\left\{\{\xi_1,\,\xi_3,\,\xi_5\},\,\{\xi_2,\,\xi_4,\,\xi\}\right\}=\Lambda^*$.
Then we have $\dot{\cal I}(t)=0$.
\end{lemma}
\noindent
{\it Proof.}
By symmetry, it suffices to consider the case $(\xi_1,\,\xi_3,\,\xi_5,\,\xi_2,\,\xi_4,\,\xi)=(\alpha_1,\,\alpha_1,\,\alpha_2,\,\beta_1,\,\beta_1,\,\beta_2)$.

First one notices that
\begin{eqnarray}\label{p-2}
& & \xi_1|d_{\xi_1}|^2-\xi_2|d_{\xi_2}|^2+\xi_3|d_{\xi_3}|^2-\xi_4|d_{\xi_4}|^2+\xi_5|d_{\xi_5}|^2-\xi|d_{\xi}|^2\nonumber\\
& = & K\left(2\alpha_1+2\beta_1+\frac{\alpha_2}{2}+\frac{\beta_2}{2}\right)-2\beta_1-\frac{\beta_2}{2}.
\end{eqnarray}
From the crucial fact that $(\alpha_1,\,\alpha_2,\,\beta_1,\,\beta_2)=(M,\,-3M-N,\,N,\,-3N-M)$, we see that the first term on the right-hand side of (\ref{p-2}) is zero.
From the identity
\begin{eqnarray}\label{t1}
\alpha_1+\beta_1=M+N,\quad 2\alpha_1+2\beta_1+\frac{\alpha_2+\beta_2}{2}=0,
\end{eqnarray}
and
\begin{eqnarray}\label{t2}
2\alpha_1-2\beta_1+\alpha_2-\beta_2=0,
\end{eqnarray}
thus $\dot{\cal I}$ simplifies to
\begin{eqnarray*}
& & 2\alpha_1^2-2\beta_1^2+\alpha_2^2-\beta_2^2+\lambda\left(2\beta_1+\frac{\beta_2}{2}\right)\\
& = & 10(\alpha_1+\beta_1)(\alpha_1-\beta_1)+\frac{\lambda}{4}(\alpha_2-\beta_2)\\
& =  & 10(\alpha_1-\beta_1)\left(\alpha_1+\beta_1-\frac{\lambda}{20}\right)=0,
\end{eqnarray*}
since $\lambda=20(\alpha_1+\beta_1)$.
This completes the proof of Lemma \ref{phase-2}.
\qed

Finally we reduce the form (\ref{d-dnls}) to the following form:
\begin{eqnarray}\label{dm-dnls}
i\dot{d}_{\xi} & = & \mu |d_{\xi}|^2(4|d_{\xi}|^2-9)d_{\xi}\nonumber\\
& & +\mu\sum_{\left\{\{\xi_1,\,\xi_3,\,\xi_5\},\,\{\xi_2,\,\xi_4,\,\xi\}\right\}=\Lambda^*} d_{\xi_1}\overline{d_{\xi_2}}d_{\xi_3}\overline{d_{\xi_4}}d_{\xi_5},
\end{eqnarray}
with initial data $d_{\xi}(0)=a_{\xi}(0)$.

\section{ODE dynamics}\label{odedynamics}
\indent
 
We now use the symplectic polar coordinates
$$
d_{\xi}(t)=\sqrt{I_{\xi}(t)}\exp(i\theta_{\xi}(t)),\quad \xi\in\Lambda,
$$
where $\theta_{\xi}(t)\in \mathbb{R}$.
From (\ref{c_3}) and (\ref{c_4}), we may assume that $0<I_{\xi}(t)<1$ (see the hypotheses in Corollary \ref{main-corollary}).
This simply reflects the constrain on %Hamiltonian 
system of (\ref{dm-dnls}).
Substituting these equations into (\ref{dm-dnls}) gives
\begin{eqnarray}\label{s-1}
-\dot{\theta}_{\xi}=\mu I_{\xi}(4I_{\xi}-9)+\mu\sum_{\left\{\{\xi_1,\,\xi_3,\,\xi_5\},\,\{\xi_2,\,\xi_4,\,\xi\}\right\}=\Lambda^*} \left(\frac{\prod_{j=1}^5I_{\xi_j}}{I_{\xi}}\right)^{1/2}\cos(\phi_{\xi}),
\end{eqnarray}
\begin{eqnarray}\label{s-2}
\frac{1}{2}\dot{I_{\xi}}=\sum_{\left\{\{\xi_1,\,\xi_3,\,\xi_5\},\,\{\xi_2,\,\xi_4,\,\xi\}\right\}=\Lambda^*}\left(\prod_{j=1}^5I_{\xi_j}\right)^{1/2}\left(I_{\xi}\right)^{1/2}\sin(\phi_{\xi}),
\end{eqnarray}
where
\begin{eqnarray*}
\phi_{\xi}=\left\{
\begin{array}{ll}
2\theta_{\alpha_1}+\theta_{\alpha_2}-2\theta_{\beta_1}-\theta_{\beta_2}, & \mbox{when $\xi\in\{\beta_1,\,\beta_2\}$},\\
2\theta_{\beta_1}+\theta_{\beta_2}-2\theta_{\alpha_1}-\theta_{\alpha_2}, & \mbox{when $\xi\in\{\alpha_1,\,\alpha_2\}$}.
\end{array}
\right.
\end{eqnarray*}
\begin{remark}\label{ham-s}
It is possible to consider the following Hamiltonian function associated to (\ref{s-1})-(\ref{s-2}): 
\begin{eqnarray}\label{ha}
H & = & -\mu\sum_{\xi\in\Lambda}\left(\frac{4}{3}I_{\xi}^3-\frac{9}{2}I_{\xi}^2\right)-6\mu I_{\alpha_1}I_{\alpha_2}^{\frac12}I_{\beta_1}I_{\beta_2}^{\frac12}\cos\left(2\theta_{\alpha_1}+\theta_{\alpha_2}-2\theta_{\beta_1}-\theta_{\beta_2}\right)\nonumber\\
& = & \frac{33}{8}\mu\left(K^2+(1-K)^2\right)+\frac32\mu K(1-K)\nonumber\\
& & -3\mu K^{\frac{3}{2}}(1-K)^{\frac{3}{2}}\cos\left(2\theta_{\alpha_1}+\theta_{\alpha_2}-2\theta_{\beta_1}-\theta_{\beta_2}\right)
\end{eqnarray}
(since by Lemma \ref{phase-2}). 
In particular, we have that the system (\ref{s-1})-(\ref{s-2}) takes the symplectic form 
\begin{eqnarray*}
\left\{
\begin{array}{l}
\dot{\theta}_{\xi}=-\frac{\partial H}{\partial I_{\xi}},\\
\dot{I}_{\xi}=\frac{\partial H}{\partial \theta_{\xi}},
\end{array}
\right.
\end{eqnarray*}
for $\xi\in\{\alpha_1,\,\alpha_2,\,\beta_1,\,\beta_2\}$.% which leads to a completely integrable system.
\end{remark}

By introducing the symplectic change of variables in a similar way to \cite{gt1}, we define
\begin{eqnarray*}
\theta={}^t(\theta_{\alpha_1},\,\theta_{\alpha_2},\,\theta_{\beta_1},\,\theta_{\beta_2}),\quad
I={}^t(I_{\alpha_1},\,I_{\alpha_2},\,I_{\beta_1},\,I_{\beta_2}),
\end{eqnarray*}
and define new variables
\begin{eqnarray*}
\phi={}^t(\phi_1,\,\phi_2,\,\phi_3,\,\phi_4),\quad 
J={}^t(J_1,\,J_2,\,J_3,\,J_4),
\end{eqnarray*}
where
$$
\phi_1=2\theta_{\alpha_1}+\theta_{\alpha_2}-2\theta_{\beta_1}-\theta_{\beta_2},
$$
$$
\phi_2=\theta_{\alpha_2},\quad \phi_3=\theta_{\beta_1},\quad\phi_4=\theta_{\beta_1},
$$
$$
J_1=\frac12I_{\alpha_1},\quad J_2=-\frac12I_{\alpha_1}+I_{\alpha_2},\quad J_3=I_{\alpha_1}+I_{\beta_1},\quad J_4=\frac12I_{\alpha_1}+I_{\beta_2}.
$$
Then the Hamiltonian flow in the action-angle coordinates (\ref{s-1})-(\ref{s-2}) satisfies
\begin{eqnarray*}
\left(
\begin{array}{c}
\phi \\
J
\end{array}
\right)
=\left(
\begin{array}{cc}
A & 0 \\
0 & {}^tA^{-1}
\end{array}
\right)
\left(
\begin{array}{c}
\theta \\
I
\end{array}
\right),
\end{eqnarray*}
where
\begin{eqnarray*}
A=\left(
\begin{array}{cccc}
2 & 1 & -2 & -1 \\
0 & 1 & 0 & 0 \\
0 & 0 & 1 & 0 \\
0 & 0 & 0 & 1 
\end{array}
\right),\quad
{}^tA^{-1}=\left(
\begin{array}{cccc}
\frac12 & 0 & 0 & 0 \\
-\frac12 & 1 & 0 & 0 \\
1 & 0 & 1 & 0 \\
\frac12 & 0 & 0 & 1 
\end{array}
\right).
\end{eqnarray*}
By (\ref{c_3})-(\ref{c_4}) and Remark \ref{ham-s}, we have
$$
J_2=0,\quad J_3=1,\quad J_4=\frac12.
$$
As discussed in Remark \ref{ham-s}, the structure of (\ref{s-1}) and (\ref{s-2}) is reminiscent of what the function $K=K(t)$ constraints dynamics of $(\phi_j,\,J_j)$ as well as $(\theta_j,\,I_j)$.
Actually, we will interest in the periodic solutions  $(\theta_j,\,I_j)$.
We achieve this by finding a periodic solution $K(t)$ for the ODE system on $(\phi_1,\,J_1)$.

For the purpose of obtaining the solutions $(\phi_1,\,J_1)$, one considers the equations for $\phi_1$
\begin{eqnarray*}%\label{e-phi}
\dot{\phi_1}& = & -2\frac{\partial H}{\partial I_{\alpha_1}}-\frac{\partial H}{\partial I_{\alpha_2}}+2\frac{\partial H}{\partial I_{\beta_1}}+\frac{\partial H}{\partial I_{\beta_2}}\nonumber\\
& = & 4\mu\left(2K^2+\left(\frac{K}{2}\right)^2-2(1-K)^2-\left(\frac{1-K}{2}\right)^2\right)\nonumber\\
& & -9\mu\left(2K+\frac{K}{2}-2(1-K)-\frac{1-K}{2}\right)\nonumber\\
& & +\mu\left(12\left(\frac{K}{2}\right)^{\frac12}(1-K)\left(\frac{1-K}{2}\right)^{\frac12}+3\frac{K(1-K)\left(\frac{1-K}{2}\right)^{\frac12}}{\left(\frac{K}{2}\right)^{\frac12}}\right.\nonumber\\
& & \left. -12K\left(\frac{K}{2}\right)^{\frac12}\left(\frac{1-K}{2}\right)^{\frac12}-3\frac{K\left(\frac{K}{2}\right)^{\frac12}(1-K)}{\left(\frac{1-K}{2}\right)^{\frac12}}\right)\nonumber\\
& = & 9\mu(1-2K)\left(\frac72+K^{\frac12}(1-K)^{\frac12}\cos \phi_1\right),
\end{eqnarray*}
and for $J_1$
\begin{eqnarray*}%\label{e-J}
\dot{J}_1 & = & \frac{\partial H}{\partial \theta_{\alpha_1}}\nonumber\\
& = & 12\mu I_{\alpha_1}I_{\alpha_2}^{\frac12}I_{\beta_1}I_{\beta_2}^{\frac12}\sin \phi_1\nonumber\\
& = & 6\mu K^{\frac32}(1-K)^{\frac32}\sin \phi_1.
\end{eqnarray*}
Since from $J_1=I_{\alpha_1}/2=K/2$, we conclude that $\phi_1,\,K$ satisfy the following system:
\begin{eqnarray}\label{e}
\left\{
\begin{array}{l}
\dot{\phi}_1=9\mu(1-2K)\left(\frac72+K^{\frac12}(1-K)^{\frac12}\cos \phi_1\right),\\
\dot{K}= 12\mu K^{\frac32}(1-K)^{\frac32}\sin \phi_1.
\end{array}
\right.
\end{eqnarray}
We will apply the argument in \cite[Proposition 4.2]{gt1} for the dynamical system (\ref{e}), and have the following proposition.

\begin{proposition}\label{periodicity}
There exists $T>0$ so that $(\phi_1,\,K)$ is a $2T$-periodic solution of (\ref{e}) and
$$
K(0)+K(T)=1.
$$
Moreover, $\mu T$ is a constant independent of $\mu$.
\end{proposition}
\noindent
{\it Proof.}
The existence of a periodic solution for (\ref{e}) is similar to the one exposed in \cite[Proposition 4.2]{gt1}.
The intuition is that in the region $(\phi_1,\, K)\in [-\pi,\,\pi]\times (0,\,1)$, the equilibrium points associated with (\ref{e}) are $\omega_0=(0,\,1/2)$ and $\omega_{\pm}=(\pm \pi,\,1/2)$.
At the point $\omega_{\pm}$, we have that the flow obeys the Hamiltonian (\ref{ha})   
$$
H=\frac{45}{16}\mu,
$$
which constraints the equation (\ref{ha}) to be
\begin{eqnarray}\label{equi}
K(1-K)\left(\frac{27}{4}+3K^{\frac12}(1-K)^{\frac12}\cos \phi_1\right)=\frac{21}{16}.
\end{eqnarray}
The equation (\ref{equi}) defines two heteroclinic orbits.
It can be concluded that the equation (\ref{equi}) has periodic solutions $(\phi_1,\,K)$. 

We use the rescaling theorem done by $t\mapsto \mu t$ to concrete $\mu T$ as a constant independent of $\mu$. 
\qed

Once one has the periodic solution $K(t)$ in Proposition \ref{periodicity}, one can obtain the following proposition.

\begin{proposition}\label{fd-dnls}
%Under the condition $\lambda=20(M+N)$, 
There exist $c_*>0$ and $c>0$ that are independent of $\lambda,\,\mu,\,M,\,N$ so that the following holds: there exist  $2T$-periodic function $K(t):\mathbb{R}\to (0,1)$ and solutions $\{c_{\xi}\}_{\xi\in\Lambda}$ to the ODE system (\ref{m-dnls})-(\ref{md-dnls}) such that $\mu T=c$.
Moreover, $K(0)\le 1/2-c_*$ and $K(T)\ge 1/2+c_*$ and for all $t\in\mathbb{R}$
$$
|c_{\alpha_1(t)}|^2=2|c_{\alpha_2}(t)|^2=K(t),\quad
|c_{\beta_1}(t)|^2=2|c_{\beta_2}(t)|^2=1-K(t).
$$
\end{proposition}

\section{Error estimates}\label{error-sec}
\indent

In this section, we shall evaluate the residual terms $a_{\xi}(t),~\xi\not\in\Lambda$ for (\ref{f-dnls}).
%The proof follows by using the bootstrap argument.

Recall $M_*=\max\{|M|,|N|\}$. 
First we show the conservation laws.
\begin{lemma}\label{cons-apri}
Consider a global solution $u(t)$ to (\ref{dnls}) as formed in (\ref{four-ansatz}) satisfying
$$
a_{\xi}(0)=0\quad \mbox{for all $\xi\not\in\Lambda$},
$$
and
$$
|a_{\alpha_1}(0)|^2=2|a_{\alpha_2}(0)|^2=K(0),\quad
|a_{\beta_1}(0)|^2=2|a_{\beta_2}(0)|^2=1-K(0).
$$
Then we have that
\begin{eqnarray}\label{error-apri}
\left\|\langle \xi\rangle a_{\xi}(t)\right\|_{L_t^{\infty}\ell_{\xi}^2}\lesssim M_*.
\end{eqnarray}
\end{lemma}
\noindent
{\it Proof.}
As we shall see in (\ref{conservation}) and $M_0=3/2$ by $L^2$-norm conservation, one can give
\begin{eqnarray}\label{a-1}
\|u(t)\|_{L^2}^2=\sum_{\xi\in\Lambda}\left(|a_{\alpha_1}(0)|^2+|a_{\alpha_2}(0)|^2+|a_{\beta_1}(0)|^2+|a_{\beta_2}(0)|^2\right)=\frac{3}{2}.
\end{eqnarray}
Also by (\ref{energy-pos}) and via the Sobolev inequality, we have
\begin{eqnarray}\label{a-2}
\|\partial_xu(t)\|_{L^2}^2\lesssim E[u](0)\lesssim  \|\partial_xu(0)\|_{L^2}^2\left(1+\|u(0)\|_{L^2}^4\right)\lesssim M_*^2.
\end{eqnarray}
The estimate (\ref{error-apri}) follows from (\ref{a-1}) and (\ref{a-2}).
\qed

On the other hand, as in Section \ref{toymodel}, setting
\begin{eqnarray}\label{g-b-a}
b_{\xi}(t)=a_{\xi}(t)e^{i\Phi(t,\xi)},
\end{eqnarray}
where
\begin{eqnarray*}
\Phi(t,\xi)& = & \int_0^t\left((\lambda\xi-6M_0\mu)|a_{\xi}(s)|^2+4\mu|a_{\xi}(s)|^4-(\lambda^2+3\mu)\sum_{\xi'\in\mathbb{Z}}|a_{\xi'}(s)|^4\right)\,ds\\
& & + 2t\left(M_0^2(\lambda^2+3\mu)+2\lambda P_0\right),
\end{eqnarray*}
leads to the equation
\begin{eqnarray}\label{ff-dnls}
i\dot{b}_{\xi} & = & % (\lambda\xi-6M_0\mu)|b_{\xi}|^2b_{\xi}+4\mu|b_{\xi}|^4b_{\xi}\nonumber\\
%& &
 -\lambda\sum_{\scriptstyle * \atop{\scriptstyle \{\xi_1,\xi_3\}\cap\{\xi_2,\xi\}=\emptyset}}\xi_2b_{\xi_1}\overline{b_{\xi_2}}b_{\xi_3}e^{i\Phi_1(t,\xi_1,\xi_2,\xi_3,\xi)}\nonumber\\
& & +\lambda^2b_{\xi}\sum_{\scriptstyle \xi_1+\xi_3=\xi_2+\xi_4 \atop{\scriptstyle \{\xi_1,\xi_3\}\cap\{\xi_2,\xi\}=\emptyset} }b_{\xi_1}\overline{b_{\xi_2}}b_{\xi_3}\overline{b_{\xi_4}}e^{i\Phi_1(t,\xi_1,\xi_2,\xi_3,\xi_4)}\nonumber\\
& & +\mu\sum_{\scriptstyle * \atop{\scriptstyle \{\xi_1,\,\xi_3,\,\xi_5\}\ne \{\xi_2,\,\xi_4,\,\xi\}}} b_{\xi_1}\overline{b_{\xi_2}}b_{\xi_3}\overline{b_{\xi_4}}b_{\xi_5}e^{i\Phi_2(t,\xi_1,\xi_2,\xi_3,\xi_4,\xi_5,\xi)},
\end{eqnarray}
where
\begin{eqnarray*}
\Phi_1=\Phi_1(t,\xi_1,\xi_2,\xi_3,\xi_4)=\sum_{j=1}^4(-1)^{j-1}\left(\xi_j^2t-\Phi(t,\xi_j)\right),
\end{eqnarray*}
%\begin{eqnarray*}
%\%Phi_2=\Phi_2(t,\xi_1,\xi_2,\xi_3,\xi_4)=2t(\xi_1-\xi_2)(\xi_1-\xi_4)-\sum_{j=1}^4(-1)^{j-1}\Phi(t,\xi_j),
%\end{eqnarray*}
and
\begin{eqnarray*}
\Phi_2=\Phi_3(t,\xi_1,\xi_2,\xi_3,\xi_4,\xi_5,\xi_6)=\sum_{j=1}^6(-1)^{j-1}\left(\xi_j^2t-\Phi(t,\xi_j)\right).
\end{eqnarray*}
Notice that by (\ref{identity}), we have for $\xi_1+\xi_3=\xi_2+\xi_4$,
$$
\sum_{j=1}^4\xi_j^2=2(\xi_1-\xi_2)(\xi_1-\xi_4).
$$

We use (\ref{a-1}) and $\ell_{\xi}^2\hookrightarrow \ell_{\xi}^1$ to see that
\begin{eqnarray*}%\label{a-1-1}
\|a_{\xi}(t)\|_{L_t^{\infty}\ell_{\xi}^1}\lesssim 1, 
\end{eqnarray*}
as well as
\begin{eqnarray}\label{a-1-2}
\|b_{\xi}(t)\|_{L_t^{\infty}\ell_{\xi}^1}\lesssim 1.
\end{eqnarray}

To consider the error estimates, we need to introduce the following definition.  
\begin{definition}
%\begin{eqnarray*}
%{\cal F}[P_Lf](\xi)=\left\{
%\begin{array}{ll}
%\widehat{f}(\xi), & \mbox{if $|\xi|\lesssim |\lambda|$},\\
%0, & \mbox{if $|\xi|\gg |\lambda|$},
%\end{array}
%\right.
%\end{eqnarray*}
%and $P_Hf=(I-P_L)f$.
For $\xi\in\mathbb{Z}$, denote $\Lambda^{**}=\Lambda^{**}(\xi)$ is the set of $\mathbb{Z}^3$ such that $(\xi_1,\xi_2,\xi_3)\not\in \Lambda^{**}$ if and only if $\xi_1+\xi_3=\xi_2+\xi$ and either the following three cases holds:
%\begin{itemize}
%\item
%at least, two of $\xi_j$ are in $\Lambda$,
%\item
%$|\xi_2-\xi_{k}|\le \frac{M_*}{100}$ holds for some $k\in\{1,3\}$.
%\end{itemize}
\begin{itemize}
\item[(i)]
$\xi_j\in\Lambda$ for all $1\le j\le 3$,
\item[(ii)]
at least, two of $\xi_j$ are not elements of $\Lambda$,
\item[(iii)]
two of $\xi_j$ are elements of $\Lambda$, one of $\xi_j$ is not elements of $\Lambda$, $|\xi-\xi_1|\ge M_*/100$ and $|\xi-\xi_3|\ge M_*/100$.
\end{itemize}
We introduce the notation:
\begin{eqnarray*}
A_L(t)=\left\|P_Lb_{\xi}(t)\right\|_{\ell_{\xi\not\in\Lambda}^1},
\end{eqnarray*}
\begin{eqnarray*}
A_H(t)=\left\|\langle\xi\rangle^{\delta}P_Hb_{\xi}(t)\right\|_{\ell_{\xi\not\in\Lambda}^1},
\end{eqnarray*}
\begin{eqnarray*}
B_L(t)=\left\|P_L\left(b_{\xi}(t)-\widetilde{b}_{\xi}(t)\right)\right\|_{\ell_{\xi\not\in\Lambda}^1},
\end{eqnarray*}
and
\begin{eqnarray*}
B_H(t)=\left\|\langle\xi\rangle^{\delta}P_H\left(b_{\xi}(t)-\widetilde{b_{\xi}}(t)\right)\right\|_{\ell_{\xi\not\in\Lambda}^1},
\end{eqnarray*}
where
$$
\widetilde{b}_{\xi}(t)=i\lambda\int_0^t\sum_{\scriptstyle *\atop{\scriptstyle \Lambda^{**}}}\xi_2b_{\xi_1}(s)\overline{b_{\xi_2}(s)}b_{\xi_3}(s)e^{2is(\xi_1-\xi_2)(\xi_1-\xi)}\,ds.
$$
\end{definition}

Now we illustrate how the error term in (\ref{ff-dnls}) would stay bounded for a reasonable amount of time.
More precisely we have the following lemma.
\begin{proposition}\label{lem-error}
Let $0<\varepsilon\le \frac{1}{M_*^{\delta}},~\delta\in [1/2,\,1)$ and $b_{\xi}(t)$ be solution to (\ref{ff-dnls}).
Assume that $M_*\gg |M+N|,~|\mu|\ll M_*^2$ and for all $|t|\le T$
\begin{eqnarray*}%\label{e-1}
A_L(t)\le \varepsilon^2,
\end{eqnarray*}
\begin{eqnarray*}%\label{e-2}
A_H(t)\le \varepsilon.
\end{eqnarray*}
Then there exists a constant $c>0$ such that for all $|t|\le T$,
\begin{eqnarray*}%\label{e-3-3}
B_L(t) \le \frac{1}{10}\varepsilon^2+c \left(\lambda^2+|\mu|\right) \varepsilon^2 t,
\end{eqnarray*}
and
\begin{eqnarray*}%\label{e-3-4}
B_H(t) \le \frac{1}{10}\varepsilon+c \left(\lambda^2+|\mu|\right) \varepsilon t.
\end{eqnarray*}
\end{proposition}

Without loss of generality, we consider only the case $t\ge 0$.
Thanks to Duhamel's  formula of (\ref{ff-dnls}), one can write the solution $a_{\xi}(t)$ of (\ref{ff-dnls}) as the integral form:
\begin{eqnarray*}%\label{duhamel}
b_{\xi}(t) = \sum_{j=1}^3N_j(t,\xi),
\end{eqnarray*}
where
%\begin{eqnarray*}
%N_1(t,\xi)  =   -i\int_0^t\left((\lambda\xi-6M_0\mu)|b_{\xi}(s)|^2b_{\xi}(s)%+2\left(M_0^2(\lambda^2+3\mu)+2\lambda P_0\right)a_{\xi}(s)\\
%& & 
%+4\mu|b_{\xi}(s)|^4b_{\xi}(s)\right)\,ds,
%-(\lambda^2+3\mu)\left(\sum_{\xi'\in\mathbb{Z}}|a_{\xi'}(s)|^4\right)a_{\xi}(s)\Big)\,ds,
%\end{eqnarray*}
\begin{eqnarray*}
N_1(t,\xi)= i\lambda\int_0^t\sum_{\scriptstyle * \atop{\scriptstyle \{\xi_1,\xi_3\}\cap\{\xi_2,\xi\}=\emptyset}}\xi_2b_{\xi_1}(s)\overline{b_{\xi_2}(s)}b_{\xi_3}(s) e^{i\Phi_1(s,\xi_1,\xi_2,\xi_3,\xi)}\,ds,
\end{eqnarray*}
\begin{eqnarray*}
N_2(t,\xi)=-i\lambda^2\int_0^t b_{\xi}(s)\sum_{\scriptstyle \xi_1+\xi_3=\xi_2+\xi_4 \atop{\scriptstyle \{\xi_1,\xi_3\}\cap\{\xi_2,\xi\}=\emptyset} }b_{\xi_1}(s)\overline{b_{\xi_2}(s)}b_{\xi_3}(s)\overline{b_{\xi_4}(s)}e^{i\Phi_1(s,\xi_1,\xi_2,\xi_3,\xi_4)}\,ds,
\end{eqnarray*}
and
\begin{eqnarray*}
N_3(t,\xi)=-i\mu\int_0^t\sum_{\scriptstyle * \atop{\scriptstyle \{\xi_1,\,\xi_3,\,\xi_5\}\ne \{\xi_2,\,\xi_4,\,\xi\}}} b_{\xi_1}(s)\overline{b_{\xi_2}(s)}b_{\xi_3}(s)\overline{b_{\xi_4}(s)}b_{\xi_5}(s)e^{i\Phi_2(s,\xi_1,\xi_2,\xi_3,\xi_4,\xi_5,\xi)}\,ds.
\end{eqnarray*}
%and $M_0=3/2$.
Moreover, if one divides the sum of $N_2$ into two parts, $(\xi_1,\xi_2,\xi_3)\in\Lambda^{**}$ and otherwise, then
\begin{eqnarray*}%\label{duhamel-1}
& & b_{\xi}(t)-i\lambda\int_0^t\sum_{\Lambda^{**}}\xi_2b_{\xi_1}(s)\overline{b_{\xi_2}(s)}b_{\xi_3}(s)e^{i\Phi_1(s,\xi_1,\xi_2,\xi_3,\xi)}\,ds\nonumber\\
& = & N_{11}(t,\xi)+N_2(t,\xi)+N_3(t,\xi),
\end{eqnarray*}
where
\begin{eqnarray*}
N_{11}(t,\xi)=i\lambda\int_0^t\sum_{(\xi_1,\xi_2,\xi_3)\not\in\Lambda^{**}}\xi_2b_{\xi_1}(s)\overline{b_{\xi_2}(s)}b_{\xi_3}(s)e^{i\Phi_1(s,\xi_1,\xi_2,\xi_3,\xi)}\,ds.
\end{eqnarray*}

%Notice that $\lambda=20(M+N)$.
%We observe that the momentum $P_0=P[u](0)$ satisfies
%\begin{eqnarray}\label{moment-c}
%|P_0| & \lesssim & \left|\sum_{\xi\in\Lambda}\xi |a_{\xi}(0)|^2\right|+|\lambda|\sum_{*}|a_{\xi_1}(0)||a_{\xi_2}(0)||a_{\xi_3}(0)||a_{\xi_4}(0)|\nonumber\\
%& \lesssim  & \left|M K(0)+(-3M-N)\frac{K(0)}{2}+N(1-K(0))+(-3N-M)\frac{1-K(0)}{2}\right|+|\lambda|\nonumber\\
%& \lesssim & |M+N|+|\lambda|\sim |\lambda|.
%\end{eqnarray}

%Let us first give the estimate of $A_1(t)$.
For this purpose, we have the several lemmas.
%\begin{lemma}\label{l-1}
%For $0\le t\le T$,
%\begin{eqnarray}\label{ea1}
%\|\langle\xi\rangle^{\delta} N_1(t,\xi)\|_{\ell_{\xi\not\in\Lambda}^1}\lesssim \left(|\lambda|+|\mu|\right)\varepsilon^3t,
%\end{eqnarray}
%for $\delta\ge 1/3$.
%\end{lemma}
%\noindent
%{\it Proof of Lemma \ref{l-1}.}
%Since $|P_0|\lesssim |\lambda|$ and $\ell_{\xi}^1\hookrightarrow \ell_{\xi}^4$, the proof can be done as follows:
%\begin{eqnarray*}
%& & \|\langle\xi\rangle^{\delta} N_1(t,\xi)\|_{\ell_{\xi\not\in\Lambda}^1}\\
%& \lesssim &
%t\left\|\langle\xi\rangle^{\delta}\left((|\lambda||\xi|+|\mu|)|b_{\xi}|^3%+(\lambda^2+|\mu|)|a_{\xi}|
%+|\mu||b_{\xi}|^5\right)\right\|_{L_T^{\infty}\ell_{\xi\not\in\Lambda}^1},
%\end{eqnarray*}
%for $0<\varepsilon<1$ and $\delta\ge 1/3$, which completes the proof of Lemma \ref{l-1}.
%\qed
Let us first consider the term $N_2$.
\begin{lemma}\label{l-3}
Assume the same hypotheses of Proposition \ref{lem-error}.
Then
\begin{eqnarray*}%\label{ea3}
\|\langle\xi\rangle^{\delta}N_2(t,\xi)\|_{\ell_{\xi\not\in\Lambda}^1}\lesssim \lambda^2 \varepsilon^2 t.
\end{eqnarray*}
\end{lemma}
\noindent
{\it Proof of Lemma \ref{l-3}.}
One employees the constrain in (\ref{res-res}).
Because $\xi_1+\xi_3=\xi_2+\xi_4$ and $\{\xi_1,\xi_3\}\cap\{\xi_2,\xi_4\}=\emptyset$, there exists some $1\le j\le 4$ such that $\xi_j\not\in\Lambda$.
Hence by (\ref{a-1-2}) we have that
\begin{eqnarray*}
\|\langle\xi\rangle^{\delta}N_2(t,\xi)\|_{\ell_{\xi\not\in\Lambda}^1}\lesssim  \lambda^2 t \|\langle\xi\rangle^{\delta}b_{\xi}\|_{L_T^{\infty}\ell_{\xi\not\in\Lambda}^1}^2\lesssim \lambda^2 \varepsilon^2 t,
\end{eqnarray*}
since
$$
\|\langle\xi\rangle^{\delta}b_{\xi}\|_{L_T^{\infty}\ell_{\xi\not\in\Lambda}^1}\le \|P_Lb_{\xi}\|_{L_T^{\infty}\ell_{\xi\not\in\Lambda}^1}+\|\langle\xi\rangle^{\delta}P_Hb_{\xi}\|_{L_T^{\infty}\ell_{\xi\not\in\Lambda}^1}\le \varepsilon^2+\varepsilon\le 2\varepsilon.
$$
\qed
\begin{remark}\label{r-1-3}
By Lemma \ref{l-3} and its proof, we obtained the following estimates:
\begin{eqnarray*}
\|P_LN_2(t,\xi)\|_{\ell_{\xi\not\in\Lambda}^1}+\left\|\langle\xi\rangle^{\delta}P_HN_2(t,\xi)\right\|_{\ell_{\xi\not\in\Lambda}^1}\lesssim \lambda^2\varepsilon^2 t.
\end{eqnarray*}
\end{remark}

Next we illustrate the estimate for $N_3$. 
Before the estimate on $N_3$, we show the following lemma.
\begin{lemma}\label{ee}
Assume that $M_*\gg |M+N|,~\xi_1+\xi_3+\xi_5=\xi_2+\xi_4+\xi_6$ and $|\xi_1^2-\xi_2^2+\xi_3^2-\xi_4^2+\xi_5^2-\xi_6^2|\ll M_*^2,~\xi_j\in\Lambda$ for all $1\le j\le 5$, then $\{\xi_1,\,\xi_3,\,\xi_5\}=\{\xi_2,\,\xi_4,\,\xi_6\}$ or $\xi_6\in\Lambda$.
\end{lemma}
\noindent
{\it Proof of Lemma \ref{ee}.}
By symmetry we assume $|M|\ge |N|$ which leads $M_*=|M|$.
Put
$$
\frac{\xi_j}{M_*}=\eta_j+o\left(\frac{M+N}{M_*}\right),
$$
for $1\le j\le 6$, where $\eta_j\in\mathbb{Z}$.
By hypothesis, $\eta_1+\eta_3+\eta_5=\eta_2+\eta_4+\eta_6,~\eta_1^2-\eta_2^2+\eta_3^2-\eta_4^2+\eta_5^2-\eta_6^2=0$.
Moreover from
$$
N=(M+N)-M,\quad -3M-N=-(M+N)-2M,\quad -3N-M=-3(M+N)+2M,
$$
we have $\eta_j\in \{1,\,-1,\,2,\,-2\}$ for $1\le j\le 5$.

When $\{\xi_1,\,\xi_3,\,\xi_5\}\cap\{\xi_2,\,\xi_4,\,\xi_6\}\not =\emptyset$, we may assume $\xi_1=\xi_2$ by symmetry.
In this case, since $\eta_3+\eta_5=\eta_4+\eta_6$, we have
$$
0=\eta_3^2+\eta_5^2-\eta_4^2-\eta_6^2=2(\eta_3-\eta_4)(\eta_3-\eta_6),
$$
which implies $\{\xi_3,\,\xi_5\}=\{\xi_4,\,\xi_6\}$.

When $\xi_6\not\in\Lambda,~\{\xi_1,\,\xi_3,\,\xi_5\}\cap\{\xi_2,\,\xi_4,\,\xi_6\}=\emptyset$ and $\{\xi_1,\,\xi_3,\,\,\xi_5\}\ne\{\xi_2,\,\xi_4,\,\xi_6\}$, then there exists pairs $\{\xi_i,\,\xi_j\}\subset  \{\xi_1,\,\xi_3,\,\xi_5\}$ or $\{\xi_i,\,\xi_j\}\subset\{\xi_2,\,\xi_4,\,\xi_6\}$ such that $i\ne j$, but $\xi_i=\xi_j$.
By symmetry we may assume $i=1$ and $j=3$.
Note
$$
2\eta_1^2+\eta_5^2=\eta_2^2+\eta_4^2+(2\eta_1+\eta_5-\eta_2-\eta_4)^2.
$$
If $\{\eta_1,\,\eta_5\}\subset\{\pm 1\}$, then $\{\eta_2,\,\eta_4\}\subset\{\pm 1\},~\eta_1=\eta_5$ and $\eta_2=\eta_4=-\eta_1$, since $3=1^2+1^2+1^2$, which implies $|\eta_6|=5$, but this contradicts with $\eta_6\in\{\pm 1\}$.
\newline
If $\eta_1\in\{\pm 1\}$ and $\eta_5\in\{\pm 2\}$, then $\eta_2=-\eta_1$ (by symmetry) and $\{\eta_4,\,\eta_6\}=\{-\eta_1,\,-\eta_5\}$, since $6=1^2+1^2+2^2$.
Then $2\eta_1+\eta_5=0$ holds, which implies $\xi_6\in\Lambda$ by (\ref{t1}) and (\ref{t2}).
This is a contradiction.
\newline
If $\eta_1\in\{\pm 2\}$ and $\eta_5\in\{\pm 1\}$, then $\eta_2=-\eta_1$ (by symmetry) and $\{\eta_4,\,\eta_6\}=\{-\eta_1,\,-\eta_5\}$ since $9$ is expressed in sum of square by $1^2+2^2+2^2$, which contradicts with $2\eta_1+\eta_5=0$ as same as above.
\newline
If $\{\eta_1,\,\eta_5\}\subset\{\pm 2\}$, then $\eta_1=\eta_5$ and $\eta_2=\eta_4=\eta_6=-\eta_1$ since $12=2^2+2^2+2^2$,  which  implies $3\eta_1= -3\eta_1$. This is a contradiction.

Hence we conclude that $\{\xi_1,\xi_3,\xi_5\}=\{\xi_2,\,\xi_4,\,\xi_6\}$.
\qed

By using this lemma, we obtain the estimate on the term $N_3$.
\begin{lemma}\label{l-4}
Assume the same hypotheses of Proposition \ref{lem-error}.
%For $0\le t\le T,~|\mu|\ll M_*^2$ and $\delta>0$,
Then,
\begin{eqnarray}\label{ea4}
\left\|\langle\xi\rangle^{\delta}P_HN_3(t,\xi)\right\|_{\ell_{\xi\not\in\Lambda}^1}\le \frac{1}{100}\varepsilon+c\left(\lambda^2+|\mu|\right)\varepsilon t.
\end{eqnarray}
\end{lemma}
\noindent
{\it Proof of Lemma \ref{l-4}.}
We split the sum into two cases:
\begin{itemize}
\item[(i)]
at least, one of $\xi_j$ are not elements of $\Lambda$,
\item[(ii)]
$\xi_j\in\Lambda$ for all $1\le j\le 5$.
\end{itemize}

In the first case (i), we see that the contribution to the left-hand side of (\ref{ea4}) is bounded by
\begin{eqnarray*}
& & c |\mu|M_*^{\delta}t\left(\|P_Lb_{\xi}\|_{L_T^{\infty}\ell_{\xi\not\in\Lambda}^1}+M_{*}^{-\delta}\left\|\langle\xi\rangle^{\delta}P_Hb_{\xi}\right\|_{L_T^{\infty}\ell_{\xi\not\in\Lambda}^1}\right)\\
& \lesssim  & |\mu|M_*^{\delta}t\left(\varepsilon^2+M_*^{-\delta}\varepsilon\right)\lesssim |\mu|\varepsilon t,
\end{eqnarray*}
where we use $\|\langle\xi\rangle^{\delta} b_{\xi}\|_{L_T^{\infty}\ell_{\xi\not\in\Lambda}^1}\lesssim M_*^{\delta}$ from (\ref{error-apri}) and (\ref{a-1}).

In the second case (ii), we may assume $\xi_j\in\Lambda$ for $1\le j\le 5$,  without loss of generality.
Observe that by Lemma \ref{ee}, if $\xi_1+\xi_3+\xi_5=\xi_2+\xi_4+\xi$ and $|\xi_1^2-\xi_2^2+\xi_3^2-\xi_4^2+\xi_5^2-\xi^2|\ll M_*^2$, then $\{\xi_1,\,\xi_3,\xi_5\}=\{\xi_2,\,\xi_4,\,\xi\}$, but it obviously contradicts with the restriction of the formula of $N_3$.
This leads to the bound
\begin{eqnarray*}
|\xi_1^2-\xi_2^2+\xi_3^2-\xi_4^2+\xi_5^2-\xi^2|\gtrsim M_*^2.
\end{eqnarray*}
The main point is that by (\ref{a-1-2}),
$$
\left|\dot{\Phi}(t,\xi)\right|\lesssim \left(|\lambda|M_*+|\mu|\right),
$$
which in particular implies that
\begin{eqnarray}\label{e4-2}
\left|\dot{\Phi}_2(t,\xi_1,\xi_2,\xi_3,\xi_4,\xi_5,\xi)\right|\gtrsim M_*^2,
\end{eqnarray}
(since $|\lambda|\ll M_*$ and $|\mu|\ll M_*^2$).
Using integration by parts along with the above inequality (\ref{e4-2}), it follows that the contribution of this case to the left-hand side of (\ref{ea4}) is bounded by
\begin{eqnarray}\label{e4-1}
& & c|\mu|\sum_* \frac{\langle\xi\rangle^{\delta}}{M_*^2}\left[|b_{\xi_1}(s)\overline{b_{\xi_2}(s)}b_{\xi_3}(s)\overline{b_{\xi_4}(s)}b_{\xi_5}(s)|\right]_{s=0}^{s=t}\nonumber\\
& & +c|\mu|\sum_{*}\frac{\langle\xi\rangle^{\delta}}{M_*^2}\int_0^t\left|\partial_s\left(b_{\xi_1}(s)\overline{b_{\xi_2}(s)}b_{\xi_3}(s)\overline{b_{\xi_4}(s)}b_{\xi_5}(s)\right)\right|%e^{2is(\xi_1^2-\xi_2^2+\xi_3^2-\xi_4^2+\xi_5^2-\xi^2)}
\,ds.%\right|
\end{eqnarray}
We use the a priori bound (\ref{a-1-2}) to obtain that the first term of (\ref{e4-1}) is bounded by
\begin{eqnarray*}
c|\mu|M_*^{-2+\delta} \|b_{\xi}\|_{L_T^{\infty}\ell_{\xi}^1}^5\lesssim |\mu|M_*^{-2+\delta}\le \frac{1}{100}\varepsilon,
\end{eqnarray*}
where $\delta<1$ and $M_*$ is large enough.
%For the second term in (\ref{e4-1}), one can use the equation (\ref{ff-dnls}).
%We divide the sum into two sub-cases:
%\begin{itemize}
%\item
%$\xi_6\not\in\Lambda$,
%\item
%
%\end{itemize}
Substituting $\partial_sb_{\xi_j}$ to (\ref{e4-1}) and using (\ref{error-apri}), we easily see that the contribution of the second term of (\ref{e4-1}) is bounded by
\begin{eqnarray*}
 c|\mu|\left(|\lambda|M_*^{-1+\delta}\|b_{\xi}\|_{L_T^{\infty}\ell_{\xi\not\in\Lambda}^1}+(\lambda^2+|\mu|)M_*^{-2+\delta}\right)t
 \lesssim  \left(\lambda^2+|\mu|\right)\varepsilon t,
\end{eqnarray*}
which completes the proof of Lemma \ref{l-4}.
\qed

Analogously, we have the following lemma, by replace $\delta=0$.
\begin{lemma}\label{r-1-2}
Assume the same hypotheses of Proposition \ref{lem-error}.
Then
\begin{eqnarray*}
\left\|P_LN_3(t,\xi)\right\|_{\ell_{\xi\not\in\Lambda}^1}\le  \frac{1}{100}\varepsilon^2+ c\left(\lambda^2+|\mu|\right)\varepsilon^2 t.
\end{eqnarray*}
\end{lemma}
\noindent
{\it Proof of Lemma \ref{r-1-2}.}
From the proof of Lemma \ref{l-4}, it follows that
$$
\left\|P_LN_3(t,\xi)\right\|_{\ell_{\xi\not\in\Lambda}^1}\le \frac{1}{100}\varepsilon^2+c |\mu|t\left(\varepsilon^2+M_*^{-\delta}\varepsilon+M_*^{-2}+|\lambda|M_*^{-1}\|b_{\xi}\|_{L_T^{\infty}\ell_{\xi\not\in\Lambda}^1}+(\lambda^2+|\mu|)M_*^{-2}\right),
$$
which is bounded by
$$
\frac{1}{100}\varepsilon^2+ c\left(\lambda^2+|\mu|\right)\varepsilon^2 t.
$$
\qed

It remains to consider the term $N_{11}$.
%Because the estimate on $N_{1}$ is obtained from $N_{11}$ in conjunction with that of $\|b_{\xi}\|_{\ell_{\xi\not\in\Lambda}^1}$, we give the proof only for the term $N_{11}$. 
\begin{lemma}\label{l-2}
Assume the same hypotheses of Proposition \ref{lem-error}.
%For $0\le t\le T,~|\mu|\ll M_*$,
Then
\begin{eqnarray}\label{ea2}
\|\langle\xi\rangle^{\delta}P_HN_{11}(t,\xi)\|_{\ell_{\xi\not\in\Lambda}^1}\le \frac{1}{100}\varepsilon+c\left(\lambda^2+|\mu|\right)\varepsilon t.
\end{eqnarray}
\end{lemma}
\noindent
{\it Proof of Lemma \ref{l-2}.}
%By using $A_2$ we have to consider the case when $(\xi_1,\xi_2,\xi_3)\not\in\Lambda^{**}$.
%Then there three cases
We divide the sum on $(\xi_1,\,\xi_2,\,\xi_3)\not\in\Lambda^{**}$ into three cases:
\begin{itemize}
\item[(i)]
$\xi_j\in\Lambda$ for all $1\le j\le 3$,
\item[(ii)]
at least, two of $\xi_j$ are not elements of $\Lambda$,
\item[(iii)]
two of $\xi_j$ are elements of $\Lambda$, one of $\xi_j$ is not elements of $\Lambda$, $|\xi-\xi_1|\ge M_*/100$ and $|\xi-\xi_3|\ge M_*/100$.
\end{itemize}

In the case (i), we see that
\begin{eqnarray}\label{gr}
|(\xi-\xi_1)(\xi-\xi_3)|=|(\xi_2-\xi_3)(\xi_2-\xi_1)|\gtrsim M_*^2,
\end{eqnarray}
since $M_*\gg |\lambda|\sim |M+N|$ and the constrain $\{\xi_1,\,\xi_3\}\cap \{\xi_2,\,\xi\}=\emptyset$.
As we saw in (\ref{e4-2}), the following inequality holds:
$$
\left|\dot{\Phi}_1(s,\xi_1,\xi_2,\xi_3,\xi)\right|\gtrsim \left|(\xi_1-\xi_2)(\xi_1-\xi)\right|.
$$
Using integration by parts, it follows that the contribution of this case to the left-hand side of (\ref{ea2}) is bounded by
\begin{eqnarray}\label{e2-1}
& & c |\lambda|\sum_{\Lambda^{**}}\left|\frac{\langle\xi\rangle^{\delta}\xi_2}{(\xi_1-\xi_2)(\xi_1-\xi)}\right|\left|\left[b_{\xi_1}(s)\overline{b_{\xi_2}(s)}b_{\xi_3}(s)\right]_{s=0}^{s=t}\right|\nonumber\\
& & +c|\lambda|\int_0^t\sum_{\Lambda^{**}}\left|\frac{\langle\xi\rangle^{\delta}\xi_2}{(\xi_1-\xi_2)(\xi_1-\xi)}\right|\left|\partial_s\left(b_{\xi_1}(s)\overline{b_{\xi_2}(s)}b_{\xi_3}(s)\right)\right|\,ds.
\end{eqnarray}
By (\ref{a-1-2}) and (\ref{gr}), we have that the first term of (\ref{e2-1}) has the bound
\begin{eqnarray}\label{e2-1-1}
c|\lambda|M_*^{-1-\delta}\le \frac{1}{100}\varepsilon,
\end{eqnarray}
where $\delta<1$ and $M_*$ is large enough.
To deal with the second term in (\ref{e2-1}), we use the equation (\ref{ff-dnls}).
By symmetry we may only consider the case that the differential operator $\partial_s$ applies to $b_{\xi_1}(s)$, and so the contribution of this case to the second term of (\ref{e2-1}) is bounded by
\begin{eqnarray}\label{e2-2}
c\sum_{j=1}^3|\lambda|\int_0^t\sum_{\Lambda^{**}}M_*^{\delta-1}\left|\dot{N}_j(s,\xi_1)\overline{b_{\xi_2}(s)}b_{\xi_3}(s)\right|\,ds.
\end{eqnarray}
Fortunately, except for $j=1$ in (\ref{e2-2}), it follows that the contribution of these cases to (\ref{e2-2}) is bounded by
\begin{eqnarray*}
& & c|\lambda|^3M_*^{\delta-1}t\varepsilon^2+c|\lambda\mu|M_*^{2\delta-1}t\left(\varepsilon^2+M_*^{-2}+|\lambda|M_*^{-1}+|\mu|M_*^{-2}\right)\\
& \lesssim & \left(\lambda^2+|\mu|\right)\varepsilon t,
\end{eqnarray*}
where $\delta<1$ and $M_*$ is large enough.
For $j=1$, $\dot{N}_1(s,\xi_1)$ is expressed as
$$
i\lambda\sum_{\scriptstyle * \atop{\scriptstyle \{\xi_1',\xi_3'\}\cap \{\xi_2',\xi_1\}=\emptyset}}\xi_2'b_{\xi_1'}(s)\overline{b_{\xi_2'}(s)}b_{\xi_3'}(s)e^{i\Phi_1(s,\xi_1',\xi_2',\xi_3',\xi_1)}.
$$
Since $\xi_1\in\Lambda$ and constrain (\ref{res-res}), we have that at least one of $\xi_j'$ is elements of $\Lambda$, thus we divide the sum into two sub-cases (using symmetry):
\begin{itemize}
\item[(i)']
two of $\xi_j'$ are not elements of $\Lambda$,
\item[(ii)']
$\xi_2'\in\Lambda,~\xi_1'\in\Lambda,~\xi_3'\not\in\Lambda$,
\item[(iii)']
$\xi_2'\not\in\Lambda$ and $\xi_j'\in\Lambda$ for $j=1,\,3$.
\end{itemize}
In the case (i)', it easily follows that the contribution to (\ref{e2-2}) is bounded by
$$
c\lambda^2M_*^{\delta}\varepsilon^2 t\lesssim \lambda^2 \varepsilon t.
$$ 
In the case (ii)', note that $\xi_1,\,\xi_1',\,\xi_2',\,\xi_3'$ satisfy
$$
\xi_1+\xi_2'=\xi_1'+\xi_3',\quad \xi_1\in\Lambda,~\xi_1'\in\Lambda,~\xi_2'\in\Lambda,~\xi_3'\not\in\Lambda,
$$
which implies that either $|\xi_3'|\lesssim |M+N|\sim |\lambda|$ or $|\xi_3'|\gtrsim M_*$ holds.
Then we have that the contribution of this case to (\ref{e2-2}) is bounded by
\begin{eqnarray*}
c\lambda^2t\left(M_*^{\delta}\|P_Lb_{\xi_3'}\|_{L_T^{\infty}\ell_{\xi_3\not\in\Lambda}^1}+\left\|\langle\xi_3'\rangle^{\delta}P_Hb_{\xi_3'}\right\|_{L_T^{\infty}\ell_{\xi_3\not\in\Lambda}^1}\right)
\lesssim \lambda^2 \varepsilon t.
\end{eqnarray*}
In the case (iii)', by
$$
|\xi_2'b_{\xi_2'}|\lesssim |\lambda|\left|P_Lb_{\xi_2'}\right|+M_*^{1-\delta}\left|\langle\xi_2'\rangle^{\delta}P_Hb_{\xi_2'}\right|,
$$
we have that the contribution to (\ref{e2-2}) is bounded by
$$
c\lambda^2 t\left(|\lambda|M_*^{\delta-1}\varepsilon^2+\varepsilon\right)\lesssim \lambda^2\varepsilon t.
$$

We now consider the case (ii).
As in the case (i)' above, it follows that the contribution of this case to (\ref{e2-2}) is bounded by
\begin{eqnarray*}
& & c|\lambda|M_*^{1+\delta}t \left(\|P_Lb_{\xi}\|_{L_T^{\infty}\ell_{\xi\not\in\Lambda}^1} +M_*^{-\delta}\left\|\langle\xi\rangle^{\delta}P_Hb_{\xi}\right\|_{L_T^{\infty}\ell_{\xi\not\in\Lambda}^1}\right)^2\\
& \lesssim & |\lambda|M_*^{1+\delta}\varepsilon^2 t\left(\varepsilon+M_*^{-\delta}\right)^2
 \lesssim  \lambda^2\varepsilon t,
\end{eqnarray*}
for $\delta\ge 1/2$.

It thus remains to establish the bound in the case (iii).
As in the case (i) above, by using (\ref{gr}), it follows that the contribution of this case to (\ref{e2-2}) has the same bound as in (\ref{e2-1}).
The estimate for the first term of (\ref{e2-2}) is same, namely we have the bound (\ref{e2-1-1}).
For the second term in (\ref{e2-2}), we divide two cases (using symmetry):
\begin{itemize}
\item[(i)'']
the differential operator $\partial_s$ falls in $b_{\xi_j}$ with $\xi_j\in\Lambda$,
\item[(ii)'']
the differential operator $\partial_s$ falls in $b_{\xi_j}$ with $\xi_j\not\in\Lambda$.
\end{itemize}
In the case (i)'', the sum in $\dot{N}_1(s,\xi_j)$ contains at least one frequency elements not in $\Lambda$.
Thus as same as in the case (ii), we have that the contribution of this case to (\ref{e2-2}) is bounded by
\begin{eqnarray*}
& & c\lambda^2M_*^{\delta}\left(\|P_Lb_{\xi}\|_{L_T^{\infty}\ell_{\xi\not\in\Lambda}^1}+M_*^{-\delta}\left\|\langle\xi\rangle^{\delta}P_Hb_{\xi}\right\|_{L_T^{\infty}\ell_{\xi\not\in\Lambda}^1}\right)^2\\
&\lesssim & \lambda^2M_*^{\delta}\varepsilon^2\left(\varepsilon+M_*^{-\delta}\right)^2
 \lesssim \lambda^2 \varepsilon t.
\end{eqnarray*}
In the case (ii)'', we also use the formula (\ref{e2-1}), but we will take the absolute value at last.
For simplicity, we may assume that
$$
\xi_1\not\in\Lambda,~\xi_2\in\Lambda,~\xi_3\in\Lambda,
$$
and the second term is represented as follows:
\begin{eqnarray}\label{e2-2-1}
c\lambda^2\left|\int_0^t\sum_{**}\frac{\langle\xi\rangle^{\delta}\xi_2\xi_2'}{(\xi-\xi_1)(\xi-\xi_3)}b_{\xi_1'}(s)\overline{b_{\xi_2'}(s)}b_{\xi_3'}(s)\overline{b_{\xi_2}(s)}b_{\xi_3}(s)e^{i\Phi_2(s,\xi_1',\xi_2',\xi_3',\xi_2,\xi_3,\xi)}\,ds\right|,
\end{eqnarray}
where $\sum_{**}$ denotes the sum over
$$
\xi_1'+\xi_3'=\xi_2'+\xi_1,\quad \xi_1+\xi_3=\xi_2+\xi,
$$
and
$$
\{\xi_1'-\xi_2'+\xi_3',\,\xi_3\}\cap \{\xi_2,\xi\}=\{\xi_1',\,\xi_3'\}\cap\{\xi_2',\,\xi_1\}=\emptyset.
$$
In the case that there exist at least one $\xi_j'\not\in\Lambda$, we will use a similar argument to the case (i) in the proof of Lemma \ref{l-4} .
Then we should consider the case that
\begin{eqnarray}\label{co}
\xi_j'\in\Lambda~ \mbox{for all $1\le j\le 3$},\quad \xi_2\in\Lambda,\quad \xi_3\in\Lambda,\quad \xi\not\in\Lambda.
\end{eqnarray}
In the case that $|\Phi_2(s,\xi_1',\xi_2',\xi_3',\xi_2,\xi_3,\xi)|\gtrsim M_*^2$, using Lemma \ref{l-4}, it follows that the contribution this case to (\ref{e2-2-1}) is bounded by
\begin{eqnarray*}
 c \lambda^2M_*^{\delta}t\left(\varepsilon^2+M_*^{-2}+|\lambda|M_*^{-1}\|b_{\xi}\|_{L_T^{\infty}\ell_{\xi\not\in\Lambda}^1}+|\mu|M_*^{-2}\right)
 \lesssim \lambda^2\varepsilon t,
\end{eqnarray*}
for $\delta<1$.
\newline
In the case that  $|\Phi_2(s,\xi_1',\xi_2',\xi_3',\xi_2,\xi_3,\xi)|\ll M_*^2$, by Lemma \ref{ee}, we see that $\{\xi_1',\,\xi_3',\,\xi_3\}=\{\xi_2',\,\xi_2,\,\xi\}$.
It follows that $\{\xi_1',\,\xi_3'\}=\{\xi_2,\,\xi\}$ and $\xi_3=\xi_2'$ (also $\Phi_2(s,\xi_1',\xi_2',\xi_3',\xi_2,\xi_3,\xi)=0$ holds), which contradicts with (\ref{co}).
%yields that the contribution of this case to (\ref{e2-2-1}) is bounded by
%\begin{eqnarray*}
%& &  c\lambda^2t \|b_{\xi}\|_{L_T^{\infty}\ell_{\xi\in\Lambda}^1}^4\left\|\langle\xi\rangle^{\delta}P_Hb_{\xi}\right\|_{L_T^{\infty}\ell_{\xi\not\in\Lambda}^1}\lesssim \lambda^2 \varepsilon t,
%\end{eqnarray*}

Hence it completes the proof of Lemma \ref{l-2}.
\qed

\begin{lemma}\label{l-2-1}
Assume the same hypotheses of Proposition \ref{lem-error}.
%For $0\le t\le T,~|\mu|\ll M_*$,
Then
\begin{eqnarray*}%\label{ea2}
\left\|P_LN_{11}(t,\xi)\right\|_{\ell_{\xi\not\in\Lambda}^1}\le \frac{1}{100}\varepsilon^2+c \lambda^2 \varepsilon^2 t.
\end{eqnarray*}
\end{lemma}
\noindent
{\it Proof of Lemma \ref{l-2-1}.}
The proof is the same as that of Lemma \ref{l-2}.
% except for only different case (ii)''.
%In the case (ii)'', we have the bound
%$$
%c\lambda^2t \|b_{\xi}\|_{L_T^{\infty}\ell_{\xi\in\Lambda}^1}^4\left\|P_Lb_{\xi}\right\|_{L_T^{\infty}\ell_{\xi\not\in\Lambda}^1}\lesssim \lambda^2 \varepsilon^2 t,
%$$
%which completes the proof of Lemma \ref{l-2-1}.
\qed
\bigskip

\noindent
{\it Proof of Proposition \ref{lem-error}.}
We are now in the position to prove Proposition \ref{lem-error}.
The proof is a consequence of Lemma \ref{l-3}, Remark \ref{r-1-3}, Lemmas \ref{l-4}, \ref{r-1-2}, \ref{l-2} and \ref{l-2-1}.
\qed

We apply the bootstrap argument to Proposition \ref{lem-error}, and have the following proposition. 
\begin{proposition}\label{prop-boot}
Let $\delta\in[1/2,]$ and $a_{\xi}(t)$ be a global solution to (\ref{f-dnls}).
Then for all $|t|\ll \min\{\frac{1}{|\lambda|M_*},\frac{1}{\lambda^2+|\mu|}\}$,
$$
M_*^{\delta}\|P_La_{\xi}(t)\|_{\ell_{\xi\not\in\Lambda}^1}+\left\|\langle\xi\rangle^{\delta}P_Ha_{\xi}(t)\right\|_{\ell_{\xi\not\in\Lambda}^1}\lesssim\frac{1}{M_*^{\delta}}.
$$
\end{proposition}
\noindent
{\it Proof.}
We shall use Proposition \ref{lem-error} with $\varepsilon=\frac{1}{M_*^{\delta}}$, which yields that for if
$$
\max\left\{M_*^{\delta}\|P_La_{\xi}(t)\|_{\ell_{\xi\not\in\Lambda}^1},\,\left\|\langle\xi\rangle^{\delta}P_Ha_{\xi}(t)\right\|_{\ell_{\xi\not\in\Lambda}^1}\right\}\le \varepsilon,
$$
for $|t|\le T$, then we have
\begin{eqnarray*}
& & \max\left\{M_*^{\delta}\|P_La_{\xi}(t)\|_{\ell_{\xi\not\in\Lambda}^1},\,\left\|\langle\xi\rangle^{\delta}P_Ha_{\xi}(t)\right\|_{\ell_{\xi\not\in\Lambda}^1}\right\}\\
& \le & \max\left\{M_*^{\delta}\left\|P_L\widetilde{b}_{\xi}(t)\right\|_{\ell_{\xi\not\in\Lambda}^1},\,\left\|P_H\langle\xi\rangle^{\delta}\widetilde{b}_{\xi}(t)\right\|_{\ell_{\xi\not\in\Lambda}^1} \right\}+\frac{\varepsilon}{10}+c\left(\lambda^2+|\mu|\right)\varepsilon t,
\end{eqnarray*}
for $|t|\le T$.
By bootstrap argument and notation (\ref{g-b-a}) combined with a continuity argument, it will suffice to show that under the same hypothesis of Proposition \ref{lem-error}, 
\begin{eqnarray}\label{c-1}
\max\left\{M_*^{\delta}\left\|P_L\widetilde{b}_{\xi}(t)\right\|_{\ell_{\xi\not\in\Lambda}^1},\,\left\|P_H\langle\xi\rangle^{\delta}\widetilde{b}_{\xi}(t)\right\|_{\ell_{\xi\not\in\Lambda}^1}\right\} \lesssim |\lambda|M_*\varepsilon t.
\end{eqnarray}
Note that $(\xi_1,\xi_2,\xi_3)\in\Lambda^{**}(\xi)$ implies that
$$
|\xi|\sim|\xi_2|\sim M_* 
$$
and there is one of $\xi_j\not\in\Lambda$ such that $(\xi_1,\xi_2,\xi_3)\in\Lambda^{**}(\xi)$, which implies that the first term of (\ref{c-1}) is eliminated.
A direct computation allows us to obtain (\ref{c-1}).
\qed

\bigskip

\noindent
{\it Proof of Theorem \ref{main-theorem}.}
The proof of the corollary follows by Propositions \ref{periodicity}, \ref{fd-dnls} and \ref{prop-boot}.
\qed

\end{document}